\documentclass[11pt]{article}
\usepackage{latexsym, amsfonts}
\newcommand{\be}{\begin{equation}}
\newcommand{\ee}{\end{equation}}
\newcommand{\bea}{\begin{eqnarray}}
\newcommand{\eea}{\end{eqnarray}}
\newcommand{\bean}{\begin{eqnarray*}}
\newcommand{\eean}{\end{eqnarray*}}
\newcommand{\brray}{\begin{array}}
\newcommand{\erray}{\end{array}}
\newcommand{\ben}{\begin{equation}{nonumber}}
\newcommand{\een}{\end{equation}{nonumber}}

\newtheorem{dfn}{Definition}[section]
\newtheorem{thm}[dfn]{Theorem}
\newtheorem{lmma}[dfn]{Lemma}
\newtheorem{ppsn}[dfn]{Proposition}
\newtheorem{crlre}[dfn]{Corollary}
\newtheorem{xmpl}[dfn]{Example}
\newtheorem{rmrk}[dfn]{Remark}
\newcommand{\bdfn}{\begin{dfn}}
\newcommand{\bthm}{\begin{thm}}
\newcommand{\blmma}{\begin{lmma}}
\newcommand{\bppsn}{\begin{ppsn}}
\newcommand{\bcrlre}{\begin{crlre}}
\newcommand{\bxmpl}{\begin{xmpl}}
\newcommand{\brmrk}{\begin{rmrk}}
\newcommand{\edfn}{\end{dfn}}
\newcommand{\ethm}{\end{thm}}
\newcommand{\elmma}{\end{lmma}}
\newcommand{\eppsn}{\end{ppsn}}
\newcommand{\ecrlre}{\end{crlre}}
\newcommand{\exmpl}{\end{xmpl}}
\newcommand{\ermrk}{\end{rmrk}}

\newcommand{\IC}{\mathbb{C}}

\newcommand{\IN}{{I\! \! N}}

\newcommand{\IR}{\mathbb{R}}

\newcommand{\cla}{{\cal A}}
\newcommand{\clb}{{\cal B}}
\newcommand{\clc}{{\cal C}}

\newcommand{\clh}{{\cal H}}

\newcommand{\clk}{{\cal K}}
\newcommand{\cll}{{\cal L}}
\newcommand{\clm}{{\cal M}}

\newcommand{\clo}{{\cal O}}

\newcommand{\clq}{{\cal Q}}

\newcommand{\cls}{{\cal S}}

\newcommand{\clu}{{\cal U}}
\newcommand{\clv}{{\cal V}}
\newcommand{\clw}{{\cal W}}

\def\a*{{\cal A}_{h,*}}
\def\B{{\cal B}(h)}
\def\B1{{\cal B}_1(h)}
\def\b{{\cal B}^{\rm s.a.}(h)}
\def\b1{{\cal B}^{\rm s.a.}_1(h)}

\newcommand{\ot}{\otimes}

\newcommand{\raro}{\rightarrow}

\def \qed {$\Box$}

\begin{document}
	\[
\]
\begin{center}
{\large {\bf Quantum isometry groups of the Podles Spheres}}\\
by\\
{\large Jyotishman Bhowmick {\footnote {The support from National Board of Higher Mathematics, India,
 is gratefully acknowledged.}} }\\
and\\
{\large Debashish Goswami {\footnote{ Partially supported by a project on `Noncommutative Geometry and Quantum Groups' funded by Indian National Science Academy.}}}\\

\end{center}

\begin{abstract}
    
  For $\mu \in (0,1), c\geq 0,$ we identify the quantum group $SO_\mu(3)$ as the universal object in the category of compact quantum groups acting 
by `orientation and volume preserving isometries' in the sense of \cite{goswami2} on the natural spectral triple on the Podles sphere $S^2_{\mu, c}$ 
constructed by Dabrowski, D'Andrea, Landi and Wagner in \cite{{Dabrowski_et_al}}. 
      
  \end{abstract} 
  
 
\section{Introduction}
In a series of articles initiated by 
\cite{goswami} and followed by \cite{jyotish}, \cite{goswami2},  we have formulated and studied  a quantum group analogue of the group of Riemannian isometries of a classical or noncommutative manifold. This was motivated by  previous work of   a number of mathematicians including Wang, Banica, Bichon and others  (see, e.g. \cite{free}, \cite{wang}, 
\cite{ban1}, \cite{ban2}, \cite{bichon}, \cite{univ1} and references therein), who have defined quantum automorphism and quantum isometry groups  of finite spaces and finite dimensional algebras. Our theory of quantum isometry groups can be viewed as a natural generalization of such quantum automorphism or isometry groups of `finite' or `discrete' structures to the continuous or smooth set-up. Clearly, such a generalization is crucial to study the quantum symmetries in noncommutative geometry, and in particular, for a good understanding of quantum group equivariant spectral triples.

   The group of Riemannian isometries of a compact Riemannian manifold $M$ can be viewed as the universal object in the
      category of all compact metrizable groups acting on $M$, with smooth and isometric action. Moreover, assume that the manifold has a spin structure (hence in particular orientable, so we can fix a choice of orientation) and $D$ denotes the conventional Dirac operator acting as an unbounded self-adjoint operator on the Hilbert space $\clh$ of square integrable spinors. Then, it can be proved that the action of a compact group $G$ on the manifold lifts as a unitary representation (possibly of some group $\tilde{G}$ which is topologically a $2$-cover of  $G$, see \cite{CD} and \cite{Dabrowski_spinors} for more details) on the Hilbert space $\clh$ which commutes with $D$ if and only if   the action on the manifold is an orientation preserving isometric action. Therefore, to define the quantum analogue of the group of orientation-preserving Riemannian isometry group of a possibly noncommutative manifold given by a spectral triple $(\cla^\infty, \clh, D)$, 
       it is reasonable to  consider a
     category  ${\bf Q}^\prime(D)$ of compact quantum groups having unitary (co-) representation, say $U$, on $\clh$,  which commutes with $D$, and  the action on $\clb(\clh)$ obtained via conjugation by $U$ maps $\cla^\infty$ into its weak closure. 
A universal object in this category, if it exists, should define the `quantum group of orientation preserving Riemannian isometries' of the underlying spectral triple. Indeed (see \cite{goswami2}), if we consider a classical spectral triple, the subcategory of the category ${\bf Q}^\prime(D)$ consisting of groups has the classical group of orientation preserving isometries as the universal object, which justifies our definition of the quantum analogue.  Unfortunately, if we consider quantum group actions, even in the finite-dimensional (but with noncommutative $\cla$) situation  the category ${\bf Q}^\prime(D)$  may often fail to have a universal object.  It turns out, however, that if we fix any suitable faithful functional $ \tau_{R} $ on $\clb(\clh)$ (to be interpreted as the choice of a `volume form') then there exists a universal object in the subcategory $ {\bf Q}^{\prime}_{R}(D) $ of ${\bf Q}^\prime(D)$ obtained by restricting the object-class to the quantum group actions which also preserve the given functional.  The  subtle point to note here is that unlike the classical group actions on $\clb(\clh)$ which always preserve the usual trace, a quantum group action may not do so. In fact, it was proved by one of the authors in \cite{goswami_rmp} that given an object $(\clq, U)$ of ${\bf Q}^\prime(D)$ (where $\clq$ is the compact quantum group and $U$ denotes its unitary co-representation on $\clh$), we can find a suitable functional $\tau_{R}$ (which typically differs from the usual trace of $\clb(\clh)$ and can have a nontrivial modularity) which  is preserved by the action of $\clq$. This makes it quite natural to work in the setting  of twisted spectral data (as defined in \cite{goswami_rmp}). It may also be mentioned that in \cite{goswami2} we have actually worked in slightly bigger category $ {\bf Q}_{R}(D) $ of so called ` quantum family of orientation and volume preserving isometries '  and deduced that the universal object in $ {\bf Q}_{R}(D) $ exists and coincides with that of $ {\bf Q}^{\prime}_{R}(D) .$

It is very important to explicitly compute the (orientation and volume preserving) quantum group of isometries for as many examples as possible. 
This programme has been successfully carried out for a number of spectral triples, including classical spheres and tori as well as their Rieffel deformations. 
The aim of the present article is to identify $ SO_{\mu} ( 3 ) $ as the quantum group of orientation and volume preserving isometries for the spectral triples on the Podles spheres $S^2_{\mu,c},$ 
 constructed by Dabrowski et al in \cite{{Dabrowski_et_al}}. Let us mention here that although the quantum groups $SO_\mu(3)$ are `deformations' of the classical $SO(3)$, these are not Rieffel deformations and so the results and techniques of \cite{goswami2} do not apply. 

Our characterization of $SO_\mu(3)$ as the quantum isometry group of a noncommutative Riemannian manifold generalizes the classical description of the group
 $SO(3)$ as the group of orientation preserving isometries of the usual Riemannian structure on the $2$-sphere. It may be mentioned here that in a very recent
 article (\cite{soltan_new}), P. M. Soltan has characterized $SO_\mu(3)$ as the universal compact quantum group acting on the finite dimensional $C^*$-algebra
 $M_2(\IC)$ such that the action preserves a functional $\omega_\mu$ defined in \cite{soltan_new}. 
 In the classical case, we have three equivalent descriptions of $SO(3)$: (a) as a quotient of $SU(2)$, (b) as the group of (orientation preserving) isometries 
 of $S^2$, and (c) as the automorphism group of $M_2$. In the quantum case the definition of $SO_\mu(3)$ is an analogue of (a), so the
 characterization of $SO_\mu(3)$  obtained in this paper  as the quantum isometry 
 group, together with Soltan's characterization, completes the generalization of all three classical descriptions of $SO(3)$.

 \section{Notations and preliminaries}
 \subsection{Basics of the theory of compact quantum groups}
 We begin by   recalling the definition of compact quantum groups and their actions from   \cite{woro}, \cite{woro1}.  A 
compact quantum group (to be abbreviated as CQG from now on)  is given by a pair $(\cls, \Delta)$, where $\cls$ is a unital separable $C^*$ algebra 
equipped
 with a unital $C^*$-homomorphism $\Delta : \cls \raro \cls \otimes \cls$ (where $\otimes$ denotes the injective tensor product)
  satisfying \\
  (ai) $(\Delta \ot id) \circ \Delta=(id \ot \Delta) \circ \Delta$ (co-associativity), and \\
  (aii) the linear spans of $ \Delta(\cls)(\cls \ot 1)$ and $\Delta(\cls)(1 \ot \cls)$ are norm-dense in $\cls \ot \cls$. \\
  It is well-known (see \cite{woro}, \cite{woro1}) that there is a canonical dense $\ast$-subalgebra $\cls_0$ of $\cls$, consisting of the matrix coefficients of
   the finite dimensional unitary (co)-representations ( to be defined below ) of $\cls$, and maps $\epsilon : \cls_0 \raro \IC$ (co-unit) and
   $\kappa : \cls_0 \raro \cls_0$ (antipode)  defined
    on $\cls_0$ which make $\cls_0$ a Hopf $\ast$-algebra. 

    A CQG $(\cls,\Delta)$ is said to (co)-act on a unital $C^*$ algebra $\clb$,
    if there is a  unital $C^*$-homomorphism (called an action) $\alpha : \clb \raro \clb \ot \cls$ satisfying the following :\\
    (bi) $(\alpha \ot id) \circ \alpha=(id \ot \Delta) \circ \alpha$, and \\
    (bii) the linear span of $\alpha(\clb)(1 \ot \cls)$ is norm-dense in $\clb \ot \cls$.\\
     \vspace{1mm}\\

  For a Hilbert $\clb$-module $E,$ (where $\clb$ is a $C^*$ algebra) we shall denote the set of adjointable $\clb$-linear maps on $E$ by $\cll(E)$. The 
  norm-closure of the linear span of the finite-rank $\clb$-linear maps on $E$, to be called the set of compact operators on $E$, will be denoted by $\clk(E)$. We note that $\cll(E)=\clm(\clk(E))$, where $\clm(\clc)$ denotes the multiplier algebra of a $C^*$-algebra $\clc$.  We shall also need the `leg-numbering' notation: for an operator $X $ in $ \clb(\clh_1 \ot \clh_2)$,   $X_{(12)}$ and $X_{(13)}$ will denote  the operators $X \ot I_{\clh_2} $ in $ \clb(\clh_1 \ot \clh_2 \ot \clh_2)$, and $\Sigma_{23} X_{12} \Sigma_{23}$ respectively, where $\Sigma_{23}$ is the unitary on $\clh_1 \ot \clh_2 \ot \clh_2$ which flips the two copies of $\clh_2$.
  
   A unitary ( co ) representation of a CQG $ ( \cls, \Delta ) $ on a Hilbert space $ \clh $ is given by a unitary  element $ U $ of  $ \clm ( \clk ( \clh ) \otimes \cls ) \equiv \cll(\clh \ot \cls)$   satisfying  $$ ({\rm  id} \otimes \Delta ) (U) = {U}_{(12)} {U}_{(13)}.$$ 
  Given a unitary representation $U$ we shall denote by $\alpha_U$ the $\ast$-homomorphism $\alpha_U(X)=U(X \ot 1_\cls){U}^*$ for $X $ belonging to $ \clb(\clh).$ We shall sometimes identify $U$ with the isometric map from the Hilbert space $\clh$ to the Hilbert module $\clh \ot \cls$ which sends a vector $\xi$ of $\clh$ to  $U(\xi \ot 1),$ and may even denote $U(\xi \ot 1)$ by $U \xi$ by a slight abuse of notation. We say that  a (possibly unbounded) operator $T$ on $\clh$ commutes with $U$ if $T \ot I$ (with the natural domain) commutes with $U$. Such an operator will also be called $U$-equivariant or $\cls$-equivariant if $U$ is understood.

 \subsection{The quantum group of orientation preserving Riemannian isometries}
 
 We briefly recall the definition of the quantum group of orientation preserving Riemannian isometries for a spectral triple (of compact type) $ ( \cla^{\infty}, \clh, D ) $ as in \cite{goswami2}.
 We consider the category ${\bf Q^\prime}(\cla^\infty, \clh, D) \equiv{\bf Q^\prime}(D)$ whose objects (to be called orientation preserving isometries) 
are the triplets $(\cls, \Delta, U)$, where    $(\cls, \Delta)$ is a CQG with a unitary representation $U$ in $\clh$, satisfying the following:\\
 (i) $U$ commutes with $D$,\\  
(ii) for every state $ \omega $ on $ \cls, $  $({\rm id} \ot \omega) \circ \alpha_U(a) $ belongs to $ ({\cla^\infty})^{\prime \prime}$ for all $a$ in $ \cla^\infty.$

The category  $ {\bf {Q}^{\prime}}(D) $ may  not have a universal object in general, as pointed out in \cite{goswami2}. In case there is a universal object, 
we shall denote it by $ \widetilde{{QISO}^+}(D)$, with the corresponding representation $U$, say, and we denote by $QISO^+(D)$ the Woronowicz 
subalgebra of $ \widetilde{{QISO}^+}(D) $ generated by the elements of the form $<\xi \ot 1, \alpha_U(a)(\eta \ot 1)>$, 
where $\xi, \eta $ belong to $ \clh, a $ belongs to $ \cla^\infty$ and $< \cdot, \cdot >$ is the $\widetilde{{QISO}^+}(D) $-valued inner product of $\clh \ot \widetilde{{QISO}^+}(D)$. The quantum group $QISO^+(D)$ will be called the quantum group of orientation-preserving Riemannian isometries of the spectral triple $(\cla^\infty, \clh, D).$ 

 Although the category ${\bf Q}^\prime(D)$ may fail to have a universal object,  we can always get a universal object in suitable subcategories which will 
be described now. Suppose that we are given an invertible positive (possibly unbounded) operator $R$ on $\clh$  which commutes with $D$.  Then we consider
 the full subcategory ${\bf Q}^\prime_R(D)$ of ${\bf Q}^\prime(D)$ by restricting the object class to those $(\cls, \Delta, U)$ for which $\alpha_U$ 
satisfies $(\tau_R \ot {\rm id})(\alpha_U(X))=\tau_R(X)1$ for all $X$ in the $\ast$-subalgebra generated by operators of the form $|\xi><\eta|$, 
where $\xi, \eta$ are eigenvectors of the operator $D$ which by assumption has discrete spectrum, and $\tau_R(X)={\rm Tr}(RX)=<\eta, R \xi>$ 
for $X=|\xi><\eta|$.  We shall call the objects of ${\bf Q}^\prime_R (D)$ orientation and ($R$-twisted) volume preserving isometries.
It is clear (see Remark 2.9 in \cite{goswami2}) that when $Re^{-tD^2}$ is trace-class for some $t>0$, the above condition is equivalent to the condition that $\alpha_U$ preserves the bounded normal functional ${\rm Tr}(\cdot ~Re^{-tD^2})$ on the whole of $\clb(\clh)$.  It is shown in \cite{goswami2} that the category ${\bf Q}^\prime_R(D)$ always admits a universal object,  to be denoted by $ \widetilde{{QISO}^+_R}(D)$, and the Woronowicz subalgebra generated  by $\{ <\xi \ot 1, \alpha_W(a) (\eta \ot 1)> : ~\xi, \eta \in \clh, a \in \cla^\infty \}$ (where $W$ is the unitary representation of $\widetilde{{QISO}^+_R}(D)$ in $\clh$) will be denoted by $QISO^+_R(D)$ and called the quantum group of orientation and ($R$-twisted) volume preserving Riemannian isometries of the spectral triple.

 \subsection{$ SU_{\mu}( 2 ) $ and the Podles spheres}
 
  Fix $\mu$ in $(0,1)$.
 The $C^*$ algebra underlying the CQG  $ SU_{\mu}( 2 ) $ is defined as the universal unital $ C^{*} $ algebra generated by $\alpha,~ \gamma $ such that 
   $ {\alpha}^{*} \alpha + {\gamma}^{*}\gamma = 1,~  \alpha {\alpha}^{*} + {\mu}^{2}{\gamma}{\gamma}^{*} = 1,~ {\gamma}{\gamma}^{*} = {\gamma}^{*}\gamma,~ \mu \gamma \alpha = \alpha \gamma,~ \mu {\gamma}^{*} \alpha = \alpha {\gamma}^{*} .$
 
 The CQG structure is given by the following fundamental  representation: $  \left ( \begin {array} {cccc}
   \alpha & - \mu {\gamma}^{*}  \\ \gamma & {\alpha}^{*} \end {array} \right ) .$
   
   The coproduct is defined by : 
   $$ \Delta ( \alpha ) = \alpha \otimes \alpha - \mu {\gamma}^{*} \otimes \gamma, $$   
   $$ \Delta ( \gamma ) = \gamma \otimes \alpha + {\alpha}^{*} \otimes \gamma. $$
   
  The Haar state of $ SU_{\mu} ( 2 ) $ will denoted by $ h $ and the corresponding G.N.S. Hilbert space will be denoted by $L^2(SU_\mu(2))$. We will   call the unital $\ast$-subalgebra of $SU_\mu(2)$ (without any norm-closure) generated by $\alpha, \gamma$ the `co-ordinate Hopf $ \ast $-algebra' of $ SU_{\mu} ( 2 )$ and denote it by $ \clo ( SU_{\mu} ( 2 )) $ as in \cite{wagner}. 
   
 \vspace{4mm}

We now recall the definition of the Podles sphere from \cite{Dabrowski_et_al} (see also the original article \cite{podles} by Podles).

For $c\geq0$, let $ t $ in $ ( 0, 1 ]$ be given by $ c = t^{-1} - t.$
Let $ [n] \equiv [ n ]_{\mu} = \frac{\mu^n - \mu^{-n}}{\mu - \mu^{-1}} , ~ n \in \IN. $

The Podles sphere  $ S^{2}_{\mu, c} $ is defined to be the universal unital $ C^{*} $ algebra generated by elements $ x_{-1}, ~ x_{0},~ x_{1} $ satisfying the relations:
\bean x_{-1} ( x_0 - t ) = \mu^2 ( x_0 - t ) x_{-1}, \eean
\bean x_1 ( x_0 - t ) = \mu^{-2} ( x_0 - t ) x_1, \eean
\bean - [ 2 ] x_{-1} x_1 + ( \mu^2 x_0 + t ) ( x_0 - t ) = {[ 2 ]}^{2} ( 1 - t ), \eean
\bean - [ 2 ] x_1 x_{-1} + ( \mu^{-2} x_0 + t ) ( x_0 - t ) = {[ 2 ]}^{2} ( 1 - t ). \eean

The involution on $ S^{2}_{\mu, c} $ is given by
\bean  x^{*}_{-1} = - {\mu}^{-1} x_1, ~~ x^{*}_{0} = x_0. \eean

We note that $ S^{2}_{\mu,c} $ as defined above is the same as $ \chi_{q, \alpha^{\prime}, \beta } $ in page 124 of \cite{klimyk} with $ q = \mu, ~ \alpha^{\prime} = t, ~ \beta = t^2 + \mu^{-2} {( \mu^2 + 1 )}^{2} ( 1 - t ) .$

Thus, from the expressions of $ x_{-1},~ x_0, ~ x_1 $ given in page 125 of \cite{klimyk}, it follows that $ S^{2}_{\mu,c} $ can be realized as a $ \ast $-subalgebra of $ SU_{\mu}( 2 ) $ by setting:
\be \label{sphere_x_-1_in_termsof_su_mu2} x_{-1} = \frac{\mu {\alpha}^{2} + \rho ( 1 + \mu^2 ) \alpha \gamma - \mu^2 \gamma^2 }{ \mu ( 1 + \mu^2 )^{\frac{1}{2}} }, \ee
\be \label{sphere_x_0_in_termsof_su_mu2} x_0 = - \mu \gamma^* \alpha + \rho ( 1 - ( 1 + \mu^2 ) \gamma^* \gamma ) - \gamma \alpha^*, \ee
\be \label{sphere_x_1_in_termsof_su_mu2} x_1 = \frac{ \mu^2 {\gamma}^{*2} - \rho \mu ( 1 + \mu^2 ) \alpha^* \gamma^* - \mu \alpha^{*2} }{( 1 + \mu^2 )^{\frac{1}{2}} }, \ee
where $ \rho^2 = \frac{\mu^2 t^2}{ ( \mu^2 + 1 )^2 ( 1 - t ) } .$

Taking
\bean A = \frac{ 1 - t^{-1} x_0}{ 1 + \mu^2}, ~~ B = \mu ( 1 + \mu^2 )^{- \frac{1}{2}} t^{-1} x_{-1}, \eean

one obtains ( see \cite{Dabrowski_et_al} ) that the $C^*$ algebra $ S^{2}_{\mu,c} $ coincides with  the original description given in \cite{podles}, i.e, the universal $ C^{*} $ algebra generated by elements $ A $ and $ B $ satisfying the relations:
$$ A^{*} = A,~  AB = \mu^{-2} BA,$$
$$ B^* B = A - A^2 + cI, ~ B B^* = \mu^2 A - \mu^4 A^2 + cI. $$
We will denote by $ \clo ( S^{2}_{\mu,c} )$ the co-ordinate $ \ast $-algebra of $ S^{2}_{\mu,c} $, i.e. the unital $\ast$-subalgebra generated by $A,B$.   

\vspace{4mm}

We recall from \cite{wagner} the  Hopf $ \ast $-algebra $ \clu_{\mu} ( su( 2 ) ) $  which is generated by elements $ F, E, K, K^{-1} $ with defining relations :
$$ K K^{-1} = K^{-1} K = 1,~ K E = \mu E K,~ F K = \mu K F, ~ E F - F E = ( \mu - \mu^{-1} )^{-1} ( K^2 - K^{-2} ) $$
with involution $ E^* = F, ~ K^* = K $ and comultiplication :
$$ \Delta ( E ) = E \otimes K + K^{-1} \otimes E,~ \Delta ( F ) = F \otimes K + K^{-1} \otimes F, ~ \Delta ( K ) = K \otimes K .$$
The counit is given by $ \epsilon ( E ) = \epsilon ( F ) = \epsilon ( K - 1 ) = 0 $ and antipode $ S ( K ) = K^{- 1}, ~ S ( E ) = - \mu E, ~ S ( F ) = - \mu^{- 1} F. $

There is a dual pairing $ \left\langle . , . \right\rangle $ of $ \clu_{\mu} ( su( 2 ) ) $ and $ \clo ( SU_{\mu} ( 2 ) ) $, for which the nonzero values of the pairing among the generators are given below:

 $ \left\langle K^{\pm 1} , \alpha^* \right\rangle = \left\langle K^{\mp 1} , \alpha \right\rangle = \mu^{\pm \frac{1}{2}}, ~  \left\langle E , \gamma \right\rangle = \left\langle F , - \mu \gamma^* \right\rangle = 1. $

The left action $ \triangleright $ and right action $ \triangleleft $ of $ \clu_{\mu} ( su( 2 ) ) $ on $ SU_{\mu} ( 2 ) $ are given by:

$ f \triangleright x =  \left\langle f , x_{(2)} \right\rangle x_{(1)},~ x \triangleleft f =  \left\langle f , x_{(1)} \right\rangle x_{(2)}, ~ x \in \clo ( SU_{\mu} ( 2 ) ), ~ f \in \clu_{\mu} ( su( 2 ) ) $ where we have used the Sweedler notation $ \Delta ( x ) = x_{( 1 )} \otimes x_{( 2 )}.$

The actions satisfy the following :

$ {( f \triangleright x )}^{*} = {S ( f )}^{*} \triangleright x^{*}, ~ {( x \triangleleft f )}^{*} = x^{*} \triangleleft {S ( f )}^{*}, ~ f \triangleright xy = ( f_{(1)} \triangleright x ) ( f_{(2)} \triangleright y ), ~ xy \triangleleft f = ( x \triangleleft f_{(1)} ) ( y \triangleleft f_{(2)} ). $

The action on generators is given by :

 $ E \triangleright \alpha = - \mu \gamma^*, ~ E  \triangleright \gamma = \alpha^*, ~ E \triangleright \gamma^* = E  \triangleright \alpha^* = 0,~ F \triangleright ( - \mu \gamma^* ) = \alpha,~ F \triangleright \alpha^* = \gamma, ~ F  \triangleright \alpha =  F \triangleright \gamma = 0, ~ K \triangleright \alpha = \mu^{- \frac{1}{2}} \alpha, ~ K \triangleright ( \gamma^* ) = \mu^{\frac{1}{2}} \gamma^*, ~ K \triangleright \gamma = \mu^{- \frac{1}{2}} \gamma,~  K \triangleright \alpha^* = \mu^{\frac{1}{2}} \alpha^* .$

$ \gamma \triangleleft E = \alpha, ~ \alpha^{*}  \triangleleft E = - \mu \gamma^{*}, ~ \alpha \triangleleft E = \gamma^{*} \triangleleft E = 0,~ \alpha \triangleleft F = \gamma,~ - \mu \gamma^* \triangleleft F = {\alpha}^{*}, ~ \gamma \triangleleft F =  \alpha^* \triangleleft F = 0, ~ \alpha \triangleleft K = \mu^{- \frac{1}{2}} \alpha, ~ \gamma^* \triangleleft K = \mu^{ - \frac{1}{2}} \gamma^{*}, ~ \gamma \triangleleft K = \mu^{ \frac{1}{2}} \gamma,~  \alpha^{*} \triangleleft K = \mu^{\frac{1}{2}} \alpha^{*} .$

\vspace{4mm}

We recall an alternative description of $ S^{2}_{\mu,c} $ from \cite{Schmudgen_wagner_crossproduct} which we are going to need.

Let
 $$ X_{c} = \mu^{\frac{1}{2}} {( \mu^{-1} - \mu )}^{-1} c^{ - \frac{1}{2}} ( 1 - K^2 ) + E K + \mu F K, ~ ~c>0,  $$
 $$X_0=1-K^2.$$
   One has  $ \Delta ( X_c ) = 1 \otimes X_c + X_c \otimes K^2 .$ 
Moreover, we have the following (\cite{Schmudgen_wagner_crossproduct}, page 9):


\bthm
 We have, 
$$ \clo ( S^{2}_{\mu, c} ) = \{ x \in \clo ( SU_{\mu} ( 2 ) ) :  x \triangleleft X_c = 0 \}. $$
 A basis of the vector space $\clo(S^2_{\mu, c})$ is given by 
$ \{ A^k, A^kB^l, A^k{B^*}^m,~k \geq 0,~ l,m > 0 \}.$
\ethm
 Thus, any element of $\clo(S^2_{\mu,c})$ can be written as a {\it finite} linear combination of elements of the form $A^k, A^k B^l, A^k{B^*}^l.$ 

Let $ \psi $ be the densely defined linear map on $ L^{2} ( SU_{\mu} ( 2 ) ) $ defined by $ \psi ( x ) = x \triangleleft X_c.$ 

\blmma

\label{sphere_subset_ker_psi}
The map $\psi$ is closable and we have $ \overline{S^{2}_{\mu,c}} \subseteq {\rm Ker}( \overline{\psi}) $ where $ \overline{\psi} $ is the closed extension of $ \psi $ and 
 $ \overline{S^{2}_{\mu,c}} $ denotes  the Hilbert subspace generated by $ S^{2}_{\mu,c} $ in $ L^{2} ( SU_{\mu} ( 2 ) ) .$ Moreover, $\clo(S^2_{\mu,c})=\clo(SU_\mu(2)) \bigcap {\rm Ker}(\overline{\psi})=\clo(SU_\mu(2)) \bigcap {\rm Ker}(\psi).$
\elmma 

{\it Proof :} From the expression of $ X_c,$ it is clear that $ \clo ( SU_{\mu} ( 2 ) ) \subseteq {\rm Dom} ( \psi^* ) $ implying that $ \psi $ is closable, hence   $ {\rm Ker}( \overline{\psi}) $ is closed. The lemma now follows from the observation that $ \clo ( S^{2}_{\mu,c} )= {\rm Ker}( \psi) \subseteq {\rm Ker}(\overline{\psi}).$  \qed

\vspace{4mm}

We end this subsection with a discussion on the CQG   $ SO_{\mu}( 3 ) $ as described in \cite{podles_subgroup}.

It is the universal unital $ C^* $ algebra generated by elements $  M,N, G, C, L $ satisfying :\\
$ L^* L = ( I - N ) ( I - \mu^{-2} N ),~ L L^* = ( I - \mu^2 N ) ( I - \mu^4 N ),~ G^* G = G G^* = N^2,~ M^{ *} M = N - N^2,~ M M^* = \mu^2 N - \mu^4 N^2,~ C^{*} C = N - N^{2},~ C C^{*} = \mu^{2} N - \mu^{4} N^{2},~ L N = \mu^{4} N L,~ G N  = N G,~ M N = \mu^2 N M,~ C N = \mu^{2} N C,~ L G = \mu^{4} G L,~ L M = \mu^2 M L,~ M G = \mu^{2} G M,~ C M = M C,~ L G^{*} = \mu^{4} G^{*} L,~ M^2 = \mu^{-1} L G,~ M^* L = \mu^{-1} ( I - N ) C,~ N^{*} = N. $

This CQG can be identified with a Woronowicz subalgebra of $SU_\mu(2)$ by taking:
$$ N = \gamma^{*} \gamma,~ M = \alpha \gamma,~ C = \alpha {\gamma}^{*},~ G = {\gamma}^{2},~ L = {\alpha}^{2} .$$
 
 The canonical action of $SU_\mu(2)$ on $S^2_{\mu,c}$, i.e. the action obtained by restricting the coproduct of $SU_\mu(2)$ to the subalgebra 
$S^2_{\mu,c}$, is actually a faithful action of $SO_\mu(3)$. On the  subspace spanned by $ \{ x_{-1},~ x_0,~ x_1 \} $ this action is given by the 
following $SO_\mu(3)$-valued $3 \times 3$-matrix: 

$Z^\prime:= \left ( \begin {array} {cccc}
  \alpha^2  & - \mu {( 1 + \mu^{- 2} )}^{\frac{1}{2}} \alpha {\gamma}^*  &  \mu^2 \gamma^{*2} \\ {( 1 + \mu^{-2} )}^{\frac{1}{2}} \alpha \gamma  & I - \mu ( \mu + \mu^{- 1} ) \gamma^* \gamma  & - \mu {( 1 + \mu^{- 2} )}^{\frac{1}{2}} \gamma^* {\alpha}^{*}   \\ \gamma^2 & {( 1 + \mu^{-2} )}^{\frac{1}{2}} \gamma \alpha^*   &  {\alpha}^{*2} \end {array} \right ) .$

\vspace{4mm}

\section{ Spectral triples on the Podles spheres and their quantum isometry groups}

\subsection{Description of the spectral triples}
We now recall the spectral triples on $ S^{2}_{\mu, c} $ discussed in \cite{Dabrowski_et_al} (see also \cite{wagner} for the case $c=0$).

Let $ s = - c^{- \frac{1}{2}} \lambda_{-}, ~ \lambda_{\pm} = \frac{1}{2} \pm {( c + \frac{1}{4} )}^{\frac{1}{2}}. $

For all $ j $ in $ \frac{1}{2} \IN,$

 $ u_j = ( {\alpha}^{*} - s \gamma^{*} ) ( \alpha^{*} - \mu^{-1} s \gamma^* )......( \alpha^* - \mu^{- 2j + 1} s \gamma^* ),$
 
 $ w_j = ( \alpha - \mu s \gamma ) ( \alpha - \mu^2 s \gamma )........( \alpha - \mu^{2j} s \gamma ), $
 
 $ u_{- j} = E^{2j} \triangleright w_j, $
 
 $ u_0 = w_0 = 1, $
 
 $ y_1 = {( 1 + \mu^{- 2} )}^{\frac{1}{2}} ( c^{\frac{1}{2}} \mu^2 \gamma^{*2} - \mu \gamma^* \alpha^* - \mu c^{\frac{1}{2}} \alpha^{*2} ), $
 
 $ N^{l}_{kj} = {\left\| F^{l - k} \triangleright ( {y_1}^{l - \left| j \right|} u_j )  \right\|}^{- 1}. $
 
 Define $ v^{l}_{k,j} = N^{l}_{k,j} F^{l - k} \triangleright ( y^{ l - \left| j \right| }_{1} u_j ), ~ l \in \frac{1}{2} \IN_{0}, ~ j,k = - l, - l + 1,......l .$
 
 Let $ \clm_{N} $ be the Hilbert subspace of $L^2(SU_\mu(2))$  with the orthonormal basis $ \{ v^{l}_{m,N} : l = \left| N \right|, ~ \left| N \right| + 1,~ ........,~ m = - l,.......l \}.$
 
 Set 
  $$ \clh = \clm_{- \frac{1}{2}} \oplus \clm_{\frac{1}{2}} ,$$ and define a 
   representation $ \pi $ of $ S^{2}_{\mu,c} $ on $ \clh $  by
 $$ \pi ( x_i ) v^{l}_{m,N} = \alpha^{-}_{i} ( l, m; N ) v^{l - 1}_{m + i, N} + \alpha^{0}_{i} ( l, m; N ) v^{l}_{m + i, N} + \alpha^{+}_{i} ( l,m; N ) v^{l + 1}_{m + i,N}, $$
where $ \alpha^{-}_{i},~ \alpha^{0}_{i}, ~ \alpha^{+}_{i} $ are as defined in \cite{Dabrowski_et_al}.

We will often identify $ \pi ( S^{2}_{\mu,c} ) $ with $ S^{2}_{\mu,c} .$

Finally by Proposition 7.2 of \cite{Dabrowski_et_al}, the following Dirac operator $ D $ gives a spectral triple $ ( \clo(S^{2}_{\mu,c}), \clh, D ) $ which we are going to work with :
$$ D ( v^{l}_{m, \pm \frac{1}{2}} ) = ( c_{1} l + c_2 ) v^{l}_{m, \mp \frac{1}{2}}  $$
where $ c_1, c_2 $ belong to $ \IR, c_1 \neq 0 .$

\vspace{4mm}

It is easy to see that the action of $SU_\mu(2)$ on itself keeps the subspace $\clh$ invariant and so induces a unitary representation, say $ U_{0} $ on $\clh.$ 

 We define a positive, unbounded operator  $ R $ on $ \clh $ by $ R ( v^{n}_{i,\pm \frac{1}{2}} ) = \mu^{- 2i} v^{n}_{i,\pm \frac{1}{2}} .$

\bppsn

\label{sphere_tau_R_=h}

$ \alpha_{U_{0}} $ preserves the $R$-twisted volume. In particular, for $x $ belonging to $ \pi(S^2_{\mu,c})$ and $t>0$, we have $h(x)=\frac{\tau_R(x)}{\tau_R(1)}$, where $\tau_R(x):={\rm Tr}(xRe^{-tD^2})$, and $h$ denotes the restriction of the Haar state of $SU_\mu(2)$ to the subalgebra
$S^2_{\mu,c}$, which is the unique $SU_\mu(2)$-invariant state on $S^2_{\mu,c}$. 
 \eppsn
 {\it Proof :} It is enough to prove that $\tau_R$ is $\alpha_{U_0}$-invariant. 
  Define  $ R_{0} ( v^{n}_{i, \pm \frac{1}{2}} )
 = \mu^{- 2i \mp 1} v^{n}_{i, \pm \frac{1}{2}},$ and note that it has been observed in \cite{goswamisuq2} that  
${\rm Tr}(R_0e^{-tD^2})<\infty $ ( for all $ t >0 $ ) and one has $$ (\tau_{R_0} \ot {\rm id})({U_0} (x \ot 1){U_0}^*)=\tau_{R_0}(x).1,$$
 for all $x $ in $ \clb(\clh)$, where $\tau_{R_0}(x)={\rm Tr}(xR_0e^{-tD^2}).$ 
 
 Let us denote by $ P_{\frac{1}{2}},P_{- \frac{1}{2}} $  the projections onto the closed subspaces generated by $ \{ v^{l}_{i, \frac{1}{2}}  \} $ and
 $ \{ v^{l}_{i, - \frac{1}{2}} \} $ respectively. Moreover, let $ \tau_{\pm }  $ be the functionals defined by $ \tau_{\pm } ( x ) = {\rm Tr} ( x R_0 P_{\pm \frac{1}{2}} e^{- t D^2} ) .$  We observe that 
  $R_0$, $e^{-tD^2}$ and $U_0$ commute with $P_{\pm \frac{1}{2}}$  so that for $x $ belonging to $ \clb(\clh),$ 
$$ (\tau_{\pm} \ot {\rm id})(\alpha_{U_0}(x))=(\tau_{R_0} \ot {\rm id})(\alpha_{U_0}(xP_{\pm \frac{1}{2}}))=\tau_{R_0}(xP_{\pm \frac{1}{2}})1=\tau_{\pm}(x)1,$$ i.e. $\tau_{\pm}$ are $\alpha_{U_0}$-invariant. 
 Moreover, since we have 
 $ R P_{\pm \frac{1}{2}} = \mu^{\pm} R_0 P_{\pm \frac{1}{2}},$ the functional $\tau_{R}$ coincides with $\mu^{- 1} \tau_+ + \mu \tau_-$, hence is $\alpha_{U_0}$-invariant.   
 \qed
 
 \bthm
 
 $ ( SU_{\mu} ( 2 ), ~ \Delta, U_0 ) $ is an object in ${\bf Q}^\prime_R(D).$
 
 \ethm
 
 {\it Proof :} The above spectral triple is  equivariant with respect to this representation (see \cite{Dabrowski_et_al}) and it preserves $\tau_R$ by  Proposition \ref{sphere_tau_R_=h}, which completes the proof. \qed

 \vspace{4mm}

We now note down some useful  facts for later use. 

\brmrk

\label{sphere_Dabrowski_eigenspace}
   
  Using the definition of $ v^{l}_{i,j} $ and $ \triangleright $, we observe :

1. The eigenspace of $ D $ corresponding to $ ( c_{1} l + c_2 ) $ and $ - ( c_{1} l + c_2 ) $ are $ {\rm span} \{ v^{l}_{m, \frac{1}{2}} + v^{l}_{m, - \frac{1}{2}} : - l \leq m \leq l \} $ and $ {\rm span} \{ v^{l}_{m, \frac{1}{2}} - v^{l}_{m, - \frac{1}{2}} : - l \leq m \leq l \} $ respectively.

2. The  eigenspace of $ \left| D \right| $ corresponding to the eigenvalue $ ( c_{1}.\frac{1}{2} + c_2 ) $ is $  {\rm span} \{ \alpha, ~ \gamma, ~ \alpha^*, ~ \gamma^* \}. $ 

  
 \ermrk 
 
 \brmrk
 
 \label{sphere_Dabrowski_pi_A,B}
 
  1. $ \pi ( A ) v^{l}_{m,N} $ belongs to $ {\rm Span}  \{ v^{l - 1}_{m,N}, ~ v^{l}_{m,N}, ~ v^{l + 1}_{m,N} \} ,$ 
  
  $ \pi ( B ) v^{l}_{m,N} $ belongs to $ {\rm Span}  \{ v^{l - 1}_{m - 1,N},~ v^{l}_{m - 1,N},~ v^{l + 1}_{m - 1,N} \},$
   
   $ \pi ( B^* ) v^{l}_{m,N} $ belongs to $ {\rm Span}  \{ v^{l - 1}_{m + 1,N},~ v^{l}_{m + 1,N},~ v^{l + 1}_{m + 1,N} \}.$
   
   \vspace{2mm}
 
 2.$ \pi ( A^k ) ( v^{l}_{m,N} ) $ belongs to $ {\rm Span}  \{ v^{l - k}_{m,N}, ~ v^{l - k + 1}_{m,N},.......,~ v^{l + k}_{m,N} \}. $
 
 \vspace{2mm}

 3.$ \pi ( A^{m^{\prime}} B^{n^{\prime}} ) ( v^{l}_{m,N} ) $ belongs to $ {\rm Span}  \{ v^{l - m^{\prime} - n^{\prime}}_{m - n^{\prime},N},~ v^{l - ( n^{\prime} + m^{\prime} - 1 )}_{m - n^{\prime},N},~.......,~ v^{l + n^{\prime} + m^{\prime}}_{m - n^{\prime}, N} \} .$

 \vspace{2mm}

 4. $ \pi ( A^{r} B^{*s} ) ( v^l_{m,N} ) $ belongs to $ {\rm Span}  \{ v^{l - s - r}_{m + s,N},~ v^{l - s - r + 1}_{m + s,N},.......v^{l + s + r}_{m + s,N} \}. $ 
 
 \ermrk

  We shall now proceed to show that $ QISO^{+}_{R}(D) $ is isomorphic with $SO_\mu(3)$. 
  Let  $ (\tilde{ \clq}, U ) $ be an object in the category ${\bf Q}^\prime_R(D)$ of CQG s acting by orientation and $R$-twisted volume preserving
 isometries on this spectral triple and $\clq$ be the Woronowicz $ C^{*} $ subalgebra of $ \tilde{\clq}$ generated by $<(\xi \ot 1), \alpha_{U}(a) (\eta \ot 1)>_{\tilde{\clq}}$, for $\xi, \eta $ in $ \clh$, $a $ in $ S^2_{\mu,c}$ (where $< \cdot, \cdot>_{\tilde{\clq}}$ is the $\tilde{\clq}$ valued inner product of $\clh \ot \tilde{\clq}$).   We shall denote $\alpha_U$ by $\phi$ from now on.  
	
	The proof has two main steps: first,  we prove that $ \phi $ is `linear', in the sense that it keeps the span of $ \{ 1, A, B, B^{*} \} $ invariant, and then we shall exploit the facts that $\phi$ is a $\ast$-homomorphism and preserves the canonical volume form on $S^2_{\mu,c}$, i.e. the restriction of the Haar state of $SU_\mu(2)$. 
\brmrk
\label{linearity_without_r}
  The first step does not make use of the fact that $\phi$ preserves the $R$-twisted volume, so linearity of the action follows for any object in the bigger category ${\bf Q^\prime}(D)$. 

\ermrk

\subsection{Linearity of the action}

For a vector $v $ in $ \clh,$ 
we shall denote by $T_v$ the map from $\clb(\clh)$   to $L^2(SU_\mu(2))$ given by $T_v(x)=xv \in \clh \subset L^2(SU_\mu(2))$.  It is clearly a continuous map with respect to the strong operator topology on $\clb(\clh)$ and the Hilbert space topology of $L^2(SU_\mu(2))$. 


  For an element $a $ in $ SU_\mu(2)$, we consider the right multiplication $R_a $ as a bounded linear map on $L^2(SU_\mu(2))$. Clearly the composition $R_a T_v$ is a continuous linear map from $\clb(\clh) $ (with the strong operator topology) to the Hilbert space $L^2(SU_\mu(2))$.  We now define  



$$ T = R_{\alpha^*}T_\alpha+\mu^2R_\gamma T_{\gamma^*} .$$ 


%






\blmma

\label{sphere_T_tilda}
For any state $\omega$ on $\tilde{\clq}$ and $x $ in $ S^2_{\mu, c}$, we have $T(\phi_\omega(x))=\phi_\omega(x) \equiv R_1(\phi_\omega(x)) $ belonging to $ \overline{S^2_{\mu,c}} \subseteq L^2(SU_\mu(2))$, where $\phi_\omega(x)=({\rm id} \ot \omega)(\phi(x)).$

\elmma

{\it Proof :} It is clear from the definition of $T$  (using $\alpha \alpha^*+\mu^2 \gamma \gamma^*=1$) that $T(x) =x \equiv R_1(x) $ for $x $ in $ S^2_{\mu,c} \subset \clb(\clh)$, where $x$ in the right hand side of the above denotes the identification of $x \in S^2_{\mu,c}$ as a vector in $L^2(SU_\mu(2))$.  Now, the lemma follows by noting that for $ x $ in $ S^{2}_{\mu,c} ,$ $\phi_\omega(x)$ belongs to  $(S^2_{\mu,c})^{\prime \prime}$, which is  the  closure of $S^2_{\mu, c}$ in the strong operator topology, and the  continuity of $T$ in this topology  discussed before. 
 \qed

\vspace{4mm}
     
     Let     
      $$ \clv^{l} = {\rm Span} \{ v^{l^{\prime}}_{i, \pm \frac{1}{2}}, -l^{\prime} \leq i \leq l^{\prime},~ l^{\prime} \leq l \}.$$ Since $ {\rm Span} \{ v^{l}_{i, \pm \frac{1}{2}}, -l \leq i \leq l \} $ is the eigenspace of $|D|$ corresponding to the eigenvalue $c_1l+c_2$, $U$ and $U^*$ must keep  $\clv^l$ invariant  for all $ l.$

 \blmma
 
 \label{sphere_Dabrowski_linearity_1}
 
  There is some finite dimensional subspace $\clv$ of $\clo(SU_\mu(2))$  such that $ {R}_{\alpha^{*}} (\phi_\omega ( A )v^{\frac{1}{2}}_{j,\pm \frac{1}{2}}),$
  
  $ {R}_{\gamma}(\phi_\omega ( A )v^{\frac{1}{2}}_{j,\pm \frac{1}{2}}) $ belong to $  \clv$ for all states $\omega$ on $\tilde{\clq}.$
 
 The same holds when $ A $ is replaced by $ B$ or $B^{*}.$ 
 
 \elmma

{\it Proof :} We prove the result for $A$ only, since a similar argument will work for $B$ and $B^*.$

We have
$ \phi ( A )( v^{\frac{1}{2}}_{j, \pm \frac{1}{2}} \otimes 1 ) = {U} ( \pi ( A ) \otimes 1 ) {{U}}^{*} ( v^{\frac{1}{2}}_{j, \pm \frac{1}{2}} \otimes 1 ) .$ 

Now,  $ {{U}}^{*} ( v^{\frac{1}{2}}_{j, \pm \frac{1}{2}} \otimes 1 ) $ belongs to $ \clv^{\frac{1}{2}} \ot \tilde{\clq}$, and then 
using the definition of $ \pi $ as well as the  Remark \ref{sphere_Dabrowski_pi_A,B}, we observe that $ ( \pi ( A ) \otimes 1 ) {{U}}^{*} ( v^{\frac{1}{2}}_{j, \pm \frac{1}{2}} \otimes 1 ) $ belongs to $ {\rm Span}\{ v^{l^{\prime}}_{j, \pm \frac{1}{2}} : - l^{\prime} \leq  j \leq l^{\prime}, l^{\prime} \leq \frac{3}{2}   \} \ot \tilde{\clq}=\clv^{\frac{3}{2}} \ot \tilde{\clq}.$ 
Again, ${U}$ keeps $\clv^{\frac{3}{2}} \ot \tilde{\clq}$ invariant, so  $R_{\alpha^*}(\phi_\omega(A)v^{\frac{1}{2}}_{\pm \frac{1}{2}}) $ belongs to $ {\rm Span}\{ v \alpha^*: v \in \clv^{\frac{3}{2}} \}$. Similarly, $R_\gamma(\phi_\omega(A)(v^{\frac{1}{2}}_{\pm \frac{1}{2}})) $ belongs to $ {\rm Span}\{ v \gamma: v \in \clv^{\frac{3}{2}} \}.$ 
So, the lemma follows for $A$ by taking $\clv ={\rm Span} \{ v\alpha^*, v \gamma :~v \in \clv^{\frac{3}{2}} \} \subset \clo(SU_\mu(2)).$ \qed

 %
%

\vspace{4mm}

Since $\alpha, \gamma^* $ belong to $ {\rm Span}\{ v^{\frac{1}{2}}_{j,\pm \frac{1}{2}} \},$ we have the following immediate corollary:
\bcrlre
\label{dab_lin}
There is a finite dimensional subspace $\clv$ of $\clo(SU_\mu(2))$ such that for every state (hence for every bounded linear functional) $\omega$ on $\tilde{\clq}$, we have $T(\phi_\omega(A)) $ belongs to $ \clv$. A similar conclusion holds for $B$ and $B^*$ as well. 
\ecrlre

\bppsn

\label{dab_lin1}

$\phi(A),~ \phi(B),~ \phi(B^*)$ belong to $\clo(S^2_{\mu,c} )\otimes_{\rm alg} \clq.$ 

\eppsn

{\it Proof :}
 We give the proof for $\phi(A)$ only, the proof for $B,B^*$ being similar. 
From the Corollary  \ref{dab_lin} and Lemma \ref{sphere_T_tilda} it follows that  for every bounded linear functional $\omega$ on $\tilde{\clq}$,  $T(\phi_\omega(A)) $ belongs to $ \clv \bigcap \overline{S^2_{\mu,c}} \subset \clo(SU_\mu(2)) \bigcap {\rm Ker}(\psi) $ and hence $ \clv \bigcap \overline{S^2_{\mu,c}} = \clv \bigcap \clo(S^2_{\mu,c})$, where $\clv$ is the finite dimensional subspace mentioned in Corollary \ref{dab_lin}.  Clearly, $\clv \bigcap \clo(S^2_{\mu,c})$ is a finite dimensional subspace of $\clo(S^2_{\mu,c})$ implying that there must be finite $m$, say, such that  for every $\omega$, $ T(\phi_\omega ( A )) $ belongs to $ {\rm Span} \{ A^k,~ A^k B^l, ~ A^k B^{*l} : 0 \leq k,l \leq m  \}.$ Denote by ${\cal W}$ the (finite dimensional) subspace of $\clb(\clh)$ spanned by $ \{ A^k,~ A^k B^l, ~ A^k B^{*l} : 0 \leq k,l \leq m \} .$  
Since for every state (and hence for every bounded linear functional) $\omega$ on $\tilde{\clq}$, we have $T(\phi_\omega(A))=R_1(\phi_\omega(A)) \equiv \phi_\omega(A)1 $, it is clear that $\phi_\omega(A) $ belongs to $ {\cal W} $ for every $\omega $ in $ {\tilde {\clq}}^*$.  Now, let us fix any faithful state  $\omega$ on the separable unital $C^*$-algebra $\tilde{\clq}$ and embed $\tilde{\clq}$ in $\clb(L^2(\clq, \omega))\equiv \clb(\clk)$. Thus, we get a canonical embedding of $\cll(\clh \ot \tilde{\clq})$ in $\clb(\clh \ot \clk)$.  Let us thus identify $\phi(A)$ as an element of $\clb(\clh \ot \clk)$, and then by choosing a countable family of elements $\{ q_1, q_2,... \}$ of $\tilde{\clq}$ which is an orthonormal basis in $\clk=L^2(\omega)$, we can write $\phi(A)$ as a  weakly convergent series of the form $\sum_{i,j=1}^\infty \phi^{ij}(A) \ot |q_i><q_j|$. But $\phi^{ij}(A)=
({\rm id} \ot \omega_{ij})(\phi(A))$, where $\omega_{ij}(\cdot)=\omega(q_i^* \cdot q_j)$. Thus, $\phi^{ij}(A) $ belongs to $ {\cal W}$ for all $i,j$, and hence the sequence $\sum_{i,j=1}^n \phi^{ij}(A) \ot |q_i><q_j| \in {\cal W} \ot \clb(\clk)$  converges weakly, and ${\cal W}$ being finite dimensional (hence weakly closed), the limit, i.e. $\phi(A)$, must belong to ${\cal W} \ot \clb(\clk)$. In other words, if $A_1,..., A_k$ denotes a basis of ${\cal W}$, we can write $\phi(A)=\sum_{i=1}^k A_i \ot B_i$ for some $B_i \in \clb(\clk).$ 

We claim that each $B_i$ must belong to $\tilde{\clq}$. For any trace-class positive operator $\rho$ in $\clh$, say of the form $\rho=\sum_j \lambda_j |e_j><e_j|$, where $\{ e_1, e_2,,...\}$ is an orthonormal basis of $\clh$ and $\lambda_j \geq 0, \sum_j \lambda_j< \infty$,  let us denote by $\psi_\rho$ the normal functional on $\clb(\clh)$ given by $x \mapsto {\rm Tr}(\rho x)$, and it is easy to see that it has a canonical extension $\tilde{\psi}_\rho:=(\psi_\rho \ot {\rm id})$ on $\cll(\clh \ot \tilde{\clq})$ given by $\tilde{\psi}_\rho(X)=\sum _j \lambda_j < e_j \ot 1, X(e_j \ot 1)>_{\tilde{\clq}}$, where $ X $ belongs to $ \cll(\clh \ot \tilde{\clq})$ and $<\cdot, \cdot, >_{\tilde{\clq}}$ denotes the $\tilde{\clq}$-valued inner product of $\clh \ot \tilde{\clq}$. Clearly, $\tilde{\psi}_\rho$ is a bounded linear map from $\cll(\clh \ot \tilde{\clq})$ to $\tilde{\clq}$. Now, since $A_1,..., A_k$ in the expression of $\phi(A)$ are linearly independent, we can choose trace class operators $ \rho_1,..., \rho_k $ such that $\psi_{\rho_i}(A_i)=1$ and $\psi_{\rho_i}(A_j)=0$ for $j \neq i$. Then, by applying $\tilde{\psi}_{\rho_i}$ on $\phi(A)$ we conclude that $B_i $ belongs to $ \tilde{\clq}$. But by definition, $\clq$ is the Woronowicz subalgebra of $\tilde{\clq}$ generated by $<\xi \ot 1, \phi(x)(\eta \ot 1)>_{\tilde{\clq}}$, with $\eta, \xi $ belonging to $ \clh$ and $x$ in $\clo(S^2_{\mu,c})$, and hence it follows that  $B_i $ belongs to $ \clq.$ 
 \qed

\bppsn

\label{sphere_Dabrowski_linearity_2}

$ \phi $ keeps the span of $ 1, A, B, B^* $ invariant.

\eppsn 

{\it Proof :} We prove the result for $ \phi ( A ) $ only, the proof for the other cases being quite similar.

Using Proposition \ref{dab_lin1}, we can write $ \phi ( A ) $ as a finite sum of the form :

 $ \sum_{k \geq 0} A^k \otimes Q_k + \sum_{m^{\prime}, n^{\prime}, n^{\prime} \neq 0 } A^{m^{\prime}} B^{n^{\prime}} \otimes R_{m^{\prime},n^{\prime}} + \sum_{r,s, s \neq 0} A^r B^{*s} \otimes R^{\prime}_{r,s} .$

Let $ \xi = v^l_{m_0,N_{0}}. $

We have that $ U ( \xi ) $ belongs to $ {\rm Span}  \{ v^{l}_{m,N}, m = - l,......l, ~ N = \pm \frac{1}{2} \} $.
 Let us write $$ {U}( \xi \otimes 1 ) = \sum_{m = -l,....l, N = \pm \frac{1}{2} } v^{l}_{m, N} \otimes q^{l}_{(m, N),(m_{0},N_{0})},$$ 
 where $q^{l}_{(m, N),(m_{0},N_{0})} $ belong to $ \clq$. 
 Since $\alpha_U$ preserves the $R$-twisted volume, we have : \be \label{rvol111} \sum_{m^{\prime},N^{\prime}} q^l_{(m,N),~(m^{\prime},N^{\prime})} q^{l*}_{(m,N),~(m^{\prime}, N^{\prime})} = 1 .\ee
  It also follows that $ U ( A \xi ) $ belongs to $ {\rm Span} \{ v^{l^{\prime}}_{m,N}, m = - l^{\prime},........l^{\prime}, l^{\prime} = l - 1, l, l + 1, ~ N = \pm \frac{1}{2} \}. $

Recalling Remark \ref{sphere_Dabrowski_pi_A,B}, we have $ \phi ( A ) {U}( \xi \otimes 1 )
 = \sum_{k, ~ m = - l,....l, N =
 \pm \frac{1}{2} } A^k v^{l}_{m,N} \otimes Q_k q^{l}_{(m,N),~(m_{0},N_{0})} + \sum_{m^{\prime},~ n^{\prime}, n^{\prime} \neq 0, ~ m = - l,....l, N = \pm \frac{1}{2} }
 A^{m^{\prime}} B^{n^{\prime}} v^{l}_{m, N} \otimes R_{m^{\prime},~ n^{\prime}}  q^{l}_{(m,N),~(m_{0},N_{0})} $ 
 
 $ + \sum_{r,s, ~ s \neq 0, ~ m = -l,...l,~ N = 
\pm \frac{1}{2}} A^r B^{*s} v^{l}_{m, N} \otimes R^{\prime}_{r,s} q^{l}_{(m,N),~(m_{0},N_{0})}.$

Let $ m^{\prime}_{0}$ denote the  largest integer $ m^{\prime}$ such that there is a nonzero coefficient of $ A^{m^{\prime}} B^{n^{\prime}}, n^\prime \geq 1 $ 
 in the expression of $ \phi ( A ).$
We claim that the coefficient of $ v^{l - m^{\prime}_{0} - n^{\prime}}_{m - n^{\prime}, N} $ in 
$ \phi ( A ) {U} ( \xi \otimes 1 )$ is $R_{m^{\prime}_{0},n^{\prime}}q^l_{(m,N),~(m_{0},N_{0})}. $

Indeed, the term $ v^{l - m^{\prime}_{0} - n^{\prime}}_{m - n^{\prime},N} $  
can arise in three ways: it can come from a term of the form $ A^{m^{\prime \prime}} B^{n^{\prime \prime}}v^l_{m,N} $ or $ A^kv^l_{m,N} $ or $ A^{r} B^{*s} v^l_{m.N}$ 
for some  $ m^{\prime \prime}, ~ n^{\prime \prime},~ k,~ r,~ s .$

In the first case, we must have $ l - m^{\prime}_{0} - n^{\prime} = l - m^{\prime \prime} - n^{\prime \prime} + t, ~ 0 \leq t \leq 2 m^{\prime \prime} $ and $ m - n^{\prime} = m - n^{\prime \prime} $ 
implying $ m ^{\prime \prime} = m^{\prime}_{0} + t $, and since  $ m^{\prime}_{0} $ is the largest integer such that $ A^{m^{\prime}_{0}} B^{n^{\prime}} $ 
appears in $ \phi ( A ) ,$ we only have the possibility $t=0$, i.e. $ v^{l - m^{\prime}_{0} - n^{\prime}}_{m - n^{\prime}, N} $ appears only in 
$ A^{m^{\prime}_{0}} B^{n^{\prime}}.$

In the second case, we have $ m - n^{\prime} = m $ implying $ n^{\prime} = 0 $ - a contradiction.
In the last case, we have $ m - n^{\prime} = m + s $ so that $ - n^{\prime} = s $ which is only possible when $ n^{\prime} = s = 0 $ which is again a contradiction.

It now follows from the above claim, using Remark \ref{sphere_Dabrowski_pi_A,B} and comparing coefficients in the equality 
$ {U} ( A \xi \otimes 1 ) =  \phi ( A ) {U} ( \xi \otimes 1 ) $, that  
 $ R_{m^{\prime}_{0}, n^{\prime}} q^{l}_{(m,N),~(m_{0},N_{0})} = 0  $ for all $ n^{\prime} \geq 1, $ for all $ m,N $ when $ m^{\prime}_{0} \geq 1 .$ Now varying $ ( m_{0}, N_{0} ) ,$ we conclude that the above holds  for all  ( $ m_{0}, N_{0} $ ).
  Using (\ref{rvol111}), we conclude that
  
   $ R_{m^{\prime}_{0}, n^{\prime}}  \sum_{m^{\prime},N^{\prime}} q^l_{(m,N), ~(m^{\prime},N^{\prime})} q^{l*}_{(m,N),~(m^{\prime},N^{\prime})} = 0 $ for all $ n^{\prime} \geq 1 ,$
   
 that is, $ R_{m^{\prime}_{0}, n^{\prime}} = 0 $  for all $ n^{\prime} \geq 1 $ 
 if $ m^{\prime}_{0} \geq 1 .$   Proceeding by induction on $ m^{\prime}_{0},$ we deduce  $ R_{m^{\prime},n^{\prime}} = 0 $ for all $ m^{\prime} \geq 1,~ n^{\prime} \geq 1 .$

Similarly, we have $ Q_k = 0  $ for all $ k \geq 2 $ and $ R^{\prime}_{r,s} = 0 $ for all $ r \geq 1,~ s \geq 1.$ 

Thus,  $ \phi ( A ) $ belongs to $ {\rm Span}  \{ 1, A, B, B^{*}, B^{2},..., B^{n},B^{*2},..., B^{*m} \}. $
But the coefficient of $ v^{l - n^{\prime}}_{m - n^{\prime},N} $ in $ \phi ( A ) {U} ( \xi \otimes 1 )$ is $ R_{0, n^{\prime}} .$ 
Arguing as before,  we conclude that $ R_{0,n^{\prime}} = 0 $ for all $ n^{\prime} \geq 2 .$ In a similar way, we can prove $ R^{\prime}_{0,n^{\prime}} = 0 $  for all $ n^{\prime} \geq 2 .$
\qed

\vspace{2mm}

 In view of the above, let us write:  
 \be \label{alphaA} \phi( A ) = 1 \otimes T_{1} + A \otimes T_{2} + B \otimes T_{3} + B^{*} \otimes T_{4}, \ee
 \be \label{alphaB}  \phi( B ) = 1 \otimes S_{1} + A \otimes S_{2} + B \otimes S_{3} + B^{*} \otimes S_{4}, \ee
 for some $ T_i, S_i $ in $ \clq.$

   \subsection{Identification of $SO_\mu(3)$ as the quantum isometry group}

In this subsection, we shall use the facts  that $ \phi $ is a $\ast$-homomorphism and it preserves the $R$-twisted volume to derive relations among $T_i, S_i$ in ( \ref{alphaA} ), ( \ref{alphaB} ).

 \blmma
 
 \label{h_invariance}
 
 $$ T_1 = \frac{1 - T_2}{1 + \mu^2}, $$ 
 $$ S_1 = \frac{- S_2}{1 + \mu^2}. $$

 \elmma

{\it Proof :} We have the expressions of $ A $ and $ B $ in terms of the $ SU_{\mu} ( 2 ) $ elements from the equations ( \ref{sphere_x_-1_in_termsof_su_mu2} ), ~( \ref{sphere_x_0_in_termsof_su_mu2} )~ and ~ ( \ref{sphere_x_1_in_termsof_su_mu2} ). From these, we note that $ h ( A ) = {( 1 + \mu^2 )}^{-1} $ and $ h ( B ) = 0.$ By recalling Proposition \ref{sphere_tau_R_=h}, we use $ ( h \otimes {\rm id} ) \phi ( A ) = h ( A ).1 $ and $ ( h \otimes {\rm id} ) \phi ( B ) = h ( B ).1 $ to have the above two equations. \qed

 \blmma 
 
 \label{spherehomomorphism A* = A}
 
$ T^{*}_{1} = T_{1},~ T^{*}_{2} = T_{2},~ T^{*}_{4} = T_{3}. $

\elmma

{\it Proof :} It  follows by comparing the coefficients of $ 1, A $ and $ B $ respectively in the equation $ \phi( A^{*} ) = \phi( A ) .$ \qed

\vspace{4mm}

\blmma

\label{spherehomomorphism B*B = A - A2}

\bean  S^{*}_{2}S_{2} + c ( 1 + \mu^2 )^{2} S^{*}_{3} S_{3} + c {( 1 + \mu^2 )}^{2} S^{*}_{4} S_4  \eean 
\be \label{spherehomomorphism B*B = A - A2 1} = ( 1 - T_2 ) ( \mu^2 + T_2 )  - c ( 1 + \mu^2 )^{2} T_3 T^{*}_3 - c {( 1 + \mu^2 )}^{2} T^{*}_{3} T_3 + c ( 1 + \mu^2 )^{2}.1, \ee

\be \label{spherehomomorphism B*B = A - A2 2} - 2 S^{*}_{2}S_{2} + ( 1 + \mu^2 ) S^{*}_{3}S_{3} + {\mu}^{2} ( 1 + \mu^2 ) S^{*}_{4}S_{4} = ( \mu^2 + 2 T_2 - 1 ) T_{2} - \mu^2 ( 1 + \mu^2 ) T_3 T^*_3 - ( 1 + \mu^2 ) T^*_3 T_3, \ee
\be \label{spherehomomorphism B*B = A - A2 3} S^{*}_{2}S_{2} - S^{*}_{3}S_{3} - {\mu}^{4} S^{*}_{4}S_{4} = - T^{2}_{2} + {\mu}^{4}T_{3}T^*_{3} + T^*_{3}T_{3}, \ee
\be \label{spherehomomorphism B*B = A - A2 5} S^{*}_{2}S_{4} + S^{*}_{3}S_{2} = - ( \mu^2 + T_2 )T^*_{3} + T^*_3 ( 1 - T_2 ), \ee
\be \label{spherehomomorphism B*B = A - A2 6} S^{*}_{2}S_{3} + {\mu}^{2}S^{*}_{4}S_{2} = - T_{2}T_{3} - {\mu}^{2}T_{3}T_{2},  \ee
\be \label{spherehomomorphism B*B = A - A2 7}  S^{*}_{4}S_{3} = - T^{2}_{3}. \ee

\elmma

{\it Proof :} It follows by comparing the coefficients of $ 1, A, A^{2}, B^{*}, AB $ and $ B^{2} $   in the equation $ \phi( B^{*} B ) = \phi( A ) - \phi( A^{2} ) + c ~ \phi ( I ) $ and then using Lemma \ref{h_invariance}   and Lemma \ref{spherehomomorphism A* = A}. \qed


\blmma

\label{spherehomomorphism BA = mu2 AB}

 \bean - S_{2}( 1 - T_{2} ) + c {( 1 + \mu^2 )}^{2} S_3 T^{*}_3 + c {( 1 + \mu^2 )}^{2} S_4 T_3 \eean 
 \be \label{spherehomomorphism BA = mu2 AB 1} =  - {\mu}^{2} ( 1 -  T_{2} ) S_{2} + c \mu^2 {( 1 + \mu^2 )}^{2} T_3 S_4 + c \mu^2 ( 1 + \mu^2 )^{2} T^{*}_{3} S_3, \ee
  
 \be \label{spherehomomorphism BA = mu2 AB A} S_{2} - 2 S_2 T_{2} + ( 1 + \mu^2 ) ( \mu^2 S_3 T^*_3 + S_4 T_3 ) = \mu^2 S_2 - 2\mu^2 T_2 S_2 + \mu^4 ( 1 + \mu^2 ) T_3 S_4 + \mu^2 ( 1 + \mu^2 ) T^*_3 S_3,  \ee 
 \be \label{spherehomomorphism BA = mu2 AB B} - S_{2}T_{3} + S_{3}( 1 - T_{2} ) = - {\mu}^{2} T_{3}S_{2} + \mu^2 ( 1 - T_2 )S_3, \ee 
 \be \label{spherehomomorphism BA = mu2 AB B*} - S_{2}T^*_{3} + S_{4}( 1 - T_{2} ) =  {\mu}^{2} ( 1 - T_2 )S_4  - \mu^2 T^*_{3}S_{2},  \ee 
 \be \label{spherehomomorphism BA = mu2 AB AB} S_{2}T_{3} + {\mu}^{2} S_{3}T_{2} = {\mu}^{2}( T_{2}S_{3} + {\mu}^{2} T_{3}S_{2} ), \ee 
 \be \label{spherehomomorphism BA = mu2 AB B2}  S_{3}T_{3} = {\mu}^{2} T_{3} S_{3}, \ee 
 \be \label{spherehomomorphism BA = mu2 AB B*2} S_{4}T^*_{3} = {\mu}^{2} T^*_{3} S_{4}. \ee
 
 \elmma
 
 {\it Proof :} It follows by equating the coefficients of $ 1, A, B, B^{*}, AB, B^{2} $ and $ {B^{*}}^{2}  $ in the equation $ \phi( BA ) = {\mu}^{2} \phi( AB ) $ and then using Lemma \ref{h_invariance}   and Lemma \ref{spherehomomorphism A* = A}. \qed

 
 \blmma
 
 \label{spherehomomorphism BB* = mu2 A - mu4 A2}
 
 
 
 
 
 \be \label{spherehomomorphism BB* = mu2 A - mu4 A2 B}  - S_{2}S^{*}_{4} - S_{3}S^{*}_{2} = {\mu}^{2} ( 1 + \mu^2 ) T_{3} - {\mu}^{4} ( 1 - T_{2} ) T_{3} - {\mu}^{4} T_{3} ( 1 - T_{2} ), \ee 
 \be \label{spherehomomorphism BB* = mu2 A - mu4 A2 AB} S_{2}S^{*}_{4} + {\mu}^{2} S_{3}S^{*}_{2} = - {\mu}^{4} T_{2}T_{3} - {\mu}^{6} T_{3}T_{2}, \ee 
 \be \label{spherehomomorphism BB* = mu2 A - mu4 A2 B2} S_{3}S^{*}_{4} = - {\mu}^{4} T^{2}_{3}. \ee

 \elmma
 
 {\it Proof :} The Lemma is proved by equating the coefficient of $ B, AB,  B^{2} $  in the equation $ \phi( BB^{*} ) = {\mu}^{2}\phi( A ) - {\mu}^{4}\phi( A^{2} ) + c ~ \phi ( I ) $ and then using Lemma \ref{h_invariance}   and Lemma \ref{spherehomomorphism A* = A}. \qed

 
 \vspace{4mm}
 
 Now, we compute the antipode, say $ \kappa $ of $ {\clq} .$

To begin with, we note that $ \{ x_{-1}, x_0, x_{1} \} $ is a set of orthogonal vectors. Moreover, they have the same norm. The first assertion being easier, we prove below the second one.

\blmma

 $ h ( x^{*}_{-1} x_{-1} ) = h ( x^{*}_{0} x_{0} ) = h ( x^{*}_{1} x_{1} ) = t^2 ( 1 - \mu^2 ) {( 1 - \mu^6 )}^{-1} [ \mu^2 + t^{-1} ( \mu^4 + 2 \mu^2 + 1 ) + t ( - \mu^4 - 2 \mu^2 - 1 ) ]. $
 
 \elmma

{\it Proof :} We have $ x^{*}_{-1} x_{-1} = t^2 \mu^{-2} ( 1 + \mu^2 ) ( A - A^2 + cI ), ~ x^{*}_{0} x_{0} = t^2 ( 1 - 2 ( 1 + \mu^2 ) A + {( 1 + \mu^2 )}^{2} A^2 ), ~ x^{*}_{1} x_{1} = t^2 ( 1 + \mu^2 ) ( \mu^2 A - \mu^4 A^2 + cI ) .$

We recall from \cite{Schmudgen_wagner_crossproduct} that  for all  bounded Borel function $ f $ on $ \sigma ( A ) ,$ 
$$ h ( f( A ) ) = \gamma_{+} \sum_{n = 0}^{\infty} f ( \lambda_{+} \mu^{2n} ) \mu^{2n} + \gamma_{-} \sum_{n = 0}^{\infty} f ( \lambda_{-} \mu^{2n} ) \mu^{2n},  $$
where $ \lambda_{+} = \frac{1}{2} + {( c + \frac{1}{4} )}^{\frac{1}{2}}, ~ \lambda_{-} = \frac{1}{2} - {( c + \frac{1}{4} )}^{\frac{1}{2}},~ \gamma_{+} = ( 1 - \mu^2 ) \lambda_{+} {( \lambda_{+} - \lambda_{-} )}^{- 1}, ~ \gamma_{-} = ( 1 - \mu^2 ) \lambda_{-} {( \lambda_{-} - \lambda_{+} )}^{- 1} .$  

The Lemma follows by applying this relation to the above expressions of $ x^{*}_{-1} x_{-1},~ x^{*}_{0} x_{0}, ~ x^{*}_{1} x_{1}.$ \qed 

\vspace{8mm}

If $ x^{\prime}_{-1},~ x^{\prime}_{0},~ x^{\prime}_{1} $ is the normalized basis corresponding to $ \{ x_{-1},~ x_0,~ x_{1} \} $, then from ( \ref{alphaA} ) and ( \ref{alphaB} ) along with the fact that each of the vectors $  x_{-1},~ x_0,~ x_{1}  $ has the same norm, it follows that

$ \phi ( x^{\prime}_{-1} ) = x^{\prime}_{-1} \otimes S_3 + x^{\prime}_{0} \otimes - {\mu}^{-1} {( 1 + {\mu}^{2} )}^{- \frac{1}{2}} S_2 + x^{\prime}_{1} \otimes - {\mu}^{-1} S_4,  $

$ \phi ( x^{\prime}_{0}  ) = x^{\prime}_{-1} \otimes - \mu {( 1 + \mu^2 )}^{\frac{1}{2}} T_3   + x^{\prime}_{0} \otimes T_2  + x^{\prime}_{1} \otimes {( 1 + \mu^2 )}^{\frac{1}{2}} T_4,  $

$ \phi ( x^{\prime}_{1} )  = x^{\prime}_{-1} \otimes - \mu S^{*}_{4} + x^{\prime}_{0} \otimes {( 1 + \mu^2 )}^{- \frac{1}{2}} S^{*}_{2}  + x^{\prime}_{1} \otimes S^{*}_{3}.  $ 

Since $ \phi $ is kept invariant by the Haar state  $ h $  of $ SU_{\mu}( 2 ) $ and $ \phi $ keeps the span of the orthonormal set $ \{ x^{\prime}_{-1},~ x^{\prime}_{0},~ x^{\prime}_{1} \} $ invariant too,  we get a unitary representation of the CQG $ \clq $ on  the span of $ \{ x^{\prime}_{-1},~ x^{\prime}_{0},~ x^{\prime}_{1} \}.$ If we denote by $Z$ the $M_3(\clq)$-valued unitary  corresponding to this unitary representation with respect to the ordered basis $ \{ x^{\prime}_{-1},~ x^{\prime}_{0},~ x^{\prime}_{1} \} $, we get by  
using $ T_4 = T^{*}_3 $ from Lemma \ref{spherehomomorphism A* = A} the following: 

$ Z = \left ( \begin {array} {cccc}
  S_{3}  & - \mu \sqrt{1 + \mu^2} T_3  & - \mu S^{*}_4 \\ \frac{ - S_{2} }{ \mu  \sqrt{1 + {\mu}^2} }  &  T_2   &   \frac{ S^{*}_{2} }{ \sqrt{1 + {\mu}^2} } \\ - {\mu}^{- 1} S_{4}   & \sqrt{1 + \mu^2} T^{*}_{3}   &  S^{*}_{3} \end {array} \right ) .$

  Recall that ( see, for example, \cite{vandaelenotes} ), the antipode $\kappa$ on the matrix elements of a finite-dimensional   unitary representation $U^\alpha \equiv ( u_{pq}^\alpha)$ is given by $\kappa (u_{pq}^\alpha ) =( u_{qp}^\alpha )^* .$ Thus,  the antipode $\kappa$ is given by :
  
 $ \kappa ( T_2 ) = T_2,~ \kappa ( T_3 ) = \frac{S^{*}_{2}}{{\mu}^{2}( 1 + \mu^2 )},~ \kappa ( S_2 ) = \mu^2 ( 1 + \mu^2 ) T^{*}_{3},~ \kappa ( S_3 ) = S^{*}_{3},~ \kappa ( S_4 ) = {\mu}^2 S_4,~ \kappa ( T^{*}_{3} ) = \frac{S_2}{1 + \mu^2},~ \kappa ( S^{*}_2 ) = ( 1 + \mu^2 ) T_3,~ \kappa ( S^{*}_{3} ) = S_3,~ \kappa ( S^{*}_{4} ) = {\mu}^{-2} S^{*}_{4}. $

Now we derive some more relations by applying the anti-homomorphism $ \kappa $ on the relations obtained earlier.

\blmma

\label{k on B*B = A - A2}



 \bean  -2 \mu^4 {( 1 + \mu^2 )}^3 T^{*}_3 T_3 + \mu^2 {( 1 + \mu^2 )}^2 S^*_3 S_3 + \mu^4 {( 1 + \mu^2 )}^2 S_4 S^*_4 \eean
  \be \label{k on B*B = A - A2 2} = \mu^2 ( 1 + \mu^2 ) T_2 ( \mu^2 + 2T_2 - 1 ) - \mu^2 S_2 S^{*}_2 - S^{*}_2 S_2, \ee
 \bean  \mu^4 {( 1 + \mu^2 )}^4 T^{*}_3 T_3 - \mu^2 {( 1 + \mu^2 )}^2 S^*_3 S_3 - \mu^6 {( 1 + \mu^2 )}^2 S_4 S^*_4 \eean
 \be \label{k on B*B = A - A2 3}  = - \mu^2 {( 1 + \mu^2 )}^2 T^2_2 + \mu^4 S_2 S^*_2 + S^*_2 S_2, \ee
\be \label{k on B*B = A - A2 5} \mu^2 {( 1 + \mu^2 )}^2 S_4 T_3 + \mu^2 {( 1 + \mu^2 )}^2 T^*_3 S_3 = - S_2 ( \mu^2 + T_2 ) + ( 1 - T_2 )S_2, \ee
\be \label{k on B*B = A - A2 7} S_4 S_3 = - \frac{- S^2_2}{\mu^2 {( 1 + \mu^2 )}^2}. \ee

\elmma

{\it Proof :} The relations follow by applying $ \kappa $ on ( \ref{spherehomomorphism B*B = A - A2 2} ), ( \ref{spherehomomorphism B*B = A - A2 3}  ), ( \ref{spherehomomorphism B*B = A - A2 5} ) and ( \ref{spherehomomorphism B*B = A - A2 7} ) respectively. \qed


\blmma

\label{k on BA = mu2 AB}

\be \label{k on BA = mu2 AB 1} - \mu^2 ( 1 - T_2 ) T^{*}_3 + c S_2 S^{*}_3 + c S^{*}_2 S_4 = - \mu^4 T^{*}_{3}( 1 - T_2 ) + c \mu^2 S_4 S^{*}_2 + c \mu^2 S^{*}_{3} S_2, \ee
\be \label{k on BA = mu2 AB B2} S_3 S_2 = \mu^2 S_2 S_3, \ee
\be \label{k on BA = mu2 AB B*2} S_2 S_4 = \mu^2 S_4 S_2, \ee
\be \label{k on BA = mu2 AB B} - S^*_2 T^*_3 + ( 1 - T_2 )S^*_3 = - \mu^2 T^*_3 S^*_2 + \mu^2 S^*_3 ( 1 - T_2 ), \ee
\be \label{k on BA = mu2 AB B*} - S_2 T^*_3 + ( 1 - T_2 ) S_4 = \mu^2 S_4 ( 1 - T_2 ) - \mu^2 T^*_3 S_2. \ee



\elmma

{\it Proof :} The relations follow by applying $ \kappa $ on (  \ref{spherehomomorphism BA = mu2 AB 1} ), ( \ref{spherehomomorphism BA = mu2 AB B2} ), ( \ref{spherehomomorphism BA = mu2 AB B*2} ), ( \ref{spherehomomorphism BA = mu2 AB B} ) and ( \ref{spherehomomorphism BA = mu2 AB B*} ) respectively. \qed


\blmma

\label{k on BB* = mu2 A - mu4 A2}



\be \label{k on BB* = mu2 A - mu4 A2 B2} S_3 S_4 = - \frac{{\mu}^2 S^2_2}{{( 1 + \mu^2 )}^2}, \ee
\be \label{k on BB* = mu2 A - mu4 A2 B} - \mu^2 {( 1 + \mu^2 )}^2 S^*_4 T^*_3 - \mu^2 {( 1 + \mu^2 )}^2 T_3 S^*_3 = \mu^2 ( 1 + \mu^2 ) S^*_2 - \mu^4 S^*_2 ( 1 - T_2 ) - \mu^4 ( 1 - T_2 )S^*_2, \ee
\be \label{k on BB* = mu2 A - mu4 A2 AB} {( 1 + \mu^2 )}^2 S^*_4 T^*_3 + \mu^2 {( 1 + \mu^2 )}^2 T_3 S^*_3 = - \mu^2 S^*_2 T_2 - \mu^4 T_2 S^*_2. \ee


\elmma

{\it Proof :} The relations follow by applying $ \kappa $ on ( \ref{spherehomomorphism BB* = mu2 A - mu4 A2 B2} ), ( \ref{spherehomomorphism BB* = mu2 A - mu4 A2 B} ) and ( \ref{spherehomomorphism BB* = mu2 A - mu4 A2 AB} ) respectively. \qed


\brmrk

\label{sphere_remark_k}

It follows from ( \ref{k on B*B = A - A2 7} ) and ( \ref{k on BB* = mu2 A - mu4 A2 B2} ) that $ \mu^4 S_4 S_3 = S_3 S_4 .$

\ermrk
























\blmma

\label{spherecheckinghomomorphism zero}  

$ S^{*}_{2} S_2 = ( 1 - T_2 ) ( \mu^2 + T_2 ). $

\elmma

{\it Proof :} Subtracting the equation obtained by multiplying $ c ( 1 + \mu^2 ) $ with ( \ref{spherehomomorphism B*B = A - A2 2} ) from ( \ref{spherehomomorphism B*B = A - A2 1} ), we have
 \bean  ( 1 + 2c ( 1 + \mu^2 ) ) S^{*}_{2} S_2 + c {( 1 + \mu^2 )}^{2} ( 1 - \mu^2 ) S^{*}_{4} S_4  \eean 
\be \label{spherecheckinghomomorphism zero_1} = ( 1 - T_2 ) ( \mu^2 + T_2 ) - c ( 1 + \mu^2 ) ( \mu^2 + 2 T_2 - 1 ) T_2 + c {( 1 + \mu^2 )}^{2} ( \mu^2 - 1 ) T_3 T^{*}_3 + c {( 1 + \mu^2 )}^{2}.1. \ee
Again, by adding ( \ref{spherehomomorphism B*B = A - A2 1} ) with $ c {( 1 + \mu^2 )}^{2} $ times ( \ref{spherehomomorphism B*B = A - A2 3} ) gives
\bean  ( 1 + c {( 1 + \mu^2 )}^{2} ) S^{*}_{2} S_2 +  c ( 1 - \mu^4 ) {( 1 + \mu^2 )}^{2} S^{*}_{4} S_4  \eean
\be \label{spherecheckinghomomorphism zero_2} =  ( 1 - T_2 ) ( \mu^2 + T_2 ) - c {( 1 + \mu^2 )}^{2} T^{2}_{2} + c {( 1 + \mu^2 )}^{2} ( {\mu}^4 - 1 ) T_3 T^{*}_3 + c {( 1 + \mu^2 )}^{2}.1. \ee
Subtracting the equation obtained by multiplying $ ( \mu^2 + 1 ) $ with ( \ref{spherecheckinghomomorphism zero_1} ) from ( \ref{spherecheckinghomomorphism zero_2} ) we obtain
\bean - ( \mu^2 + c {( 1 + \mu^2 )}^{2} ) S^{*}_{2} S_2   
 = ( 1 - T_2 ) ( \mu^2 + T_2 ) - c {( 1 + \mu^2 )}^{2} T^{2}_{2} \eean 
 \bean  - ( 1 + \mu^2 ) ( 1 - T_2 ) ( \mu^2 + T_2 ) - c \mu^2 {( 1 + \mu^2 )}^{2}.1 + c {( 1 + \mu^2 )}^{2} ( \mu^2 + 2 T_2 - 1 ) T_2. \eean
The right hand side can be seen to equal $ - ( \mu^2 + c {( 1 + \mu^2 )}^{2} ) ( 1 - T_2 ) ( \mu^2 + T_2 ) .$

Thus, $ S^{*}_{2} S_2 = ( 1 - T_2 ) ( \mu^2 + T_2 ) .$ \qed

\blmma

\label{spherecheckinghomomorphism one}

\be \label{spherecheckinghomomorphism one_1} \mu^2 {( 1 + \mu^2 )}^{2} T^{*}_{3} T_3 = ( 1 - T_2 ) ( \mu^2 + T_2 ), \ee
\be \label{spherecheckinghomomorphism one_2} {( 1 + \mu^2 )}^{2} T_3 T^{*}_{3} = ( 1 - T_2 ) ( 1 + \mu^2 T_2 ), \ee
\be \label{spherecheckinghomomorphism one_3} S_2 S^{*}_{2} = \mu^2 ( 1 - T_2 ) ( 1 + \mu^2 T_2 ). \ee

\elmma

{\it Proof :} Applying $ \kappa $ on  Lemma \ref{spherecheckinghomomorphism zero}, we obtain ( \ref{spherecheckinghomomorphism one_1}  ).

Unitarity of the matrix $ Z $ ( ( 2, 2 ) position of the matrix $ Z^* Z $ ) gives $ \mu^2 ( 1 + \mu^2 ) T^{*}_{3} T_{3} + T^{2}_{2} + ( 1 + \mu^2 ) T_3 T^{*}_{3} = 1 .$

Using ( \ref{spherecheckinghomomorphism one_1} ) we deduce $ - ( 1 + \mu^2 )^{2} T_3 T^{*}_{3} = ( T_2 - 1 ) ( 1 + \mu^2 T_2 ).$ Thus we obtain ( \ref{spherecheckinghomomorphism one_2} ).

Applying $ \kappa $ on ( \ref{spherecheckinghomomorphism one_2} ), we deduce ( \ref{spherecheckinghomomorphism one_3}  ).  \qed

\blmma

\label{spherecheckinghomomorphism two}

$ S^{*}_{4} S_4 = S_4 S^{*}_{4} = {( 1 + \mu^2 )}^{-2} \mu^{2} {( 1 - T_2 )}^{2}. $

\elmma

{\it Proof :} Adding ( \ref{k on B*B = A - A2 2} ) and ( \ref{k on B*B = A - A2 3} ), we have :
$$ - \mu^4 {( 1 + \mu^2 )}^3 ( 1 - \mu^2 ) T^*_3 T_3 + \mu^4 {( 1 + \mu^2 )}^2 ( 1 - \mu^2 ) S_4 S^*_4 = - \mu^2 ( 1 + \mu^2 ) ( 1 - \mu^2 ) T_2 ( 1 - T_2 ) - \mu^2 ( 1 - \mu^2 ) S_2 S^*_2 .$$
Using $ \mu^2 \neq 1  ,$ we obtain,
$$ - \mu^4 {( 1 + \mu^2 )}^3 T^*_3 T_3 + \mu^4 {( 1 + \mu^2 )}^2 S_4 S^*_4 = - \mu^2 ( 1 + \mu^2 ) T_2 ( 1 - T_2 ) - \mu^2 S_2 S^*_2 .  $$
Now using ( \ref{spherecheckinghomomorphism one_1} ) and ( \ref{spherecheckinghomomorphism one_3} ), we reduce the above equation to
 \bean \lefteqn{ \mu^4 {( 1 + \mu^2 )}^2 S_4 S^*_4 }\\
 & = & - \mu^2 ( 1 - T_2 ) ( T_2 + \mu^2 T_2 + \mu^2 + \mu^4 T_2 ) + \mu^2 ( 1 + \mu^2 ) ( 1 - T_2 ) ( \mu^2 + T_2 ) \\ & = & \mu^6 {( 1 - T_2 )}^2. \eean
Thus,\bean \lefteqn{ S_4 S^*_4 }\\
& = & \frac{\mu^6}{\mu^4 {( 1 + \mu^2 )}^2} {( 1 - T_2 )}^2 \\
& = & \frac{\mu^2}{{( 1 + \mu^2 )}^2} ( 1 - T_2 )^2. \eean
Applying $ \kappa ,$ we have $ S^*_4 S_4 = \frac{\mu^2}{{( 1 + \mu^2 )}^2} {( 1 - T_2 )}^2 .$

Thus, $ S^*_4 S_4 = S_4 S^*_4 = \frac{\mu^2}{{( 1 + \mu^2 )}^2} ( 1 - T_2 )^2 .$ \qed

 




\blmma

\label{spherecheckinghomomorphism three}

$ \mu^{2} {( 1 + \mu^2 )}^{2} S^{*}_{3} S_3 = ( \mu^2 + T_2 ) [ \mu^2 ( 1 + \mu^2 ) - ( 1 - T_2 ) ]. $ 

\elmma

{\it Proof :} Using Lemma \ref{spherecheckinghomomorphism zero} in ( \ref{spherehomomorphism B*B = A - A2 1}  ), we have 
 \be \label{spherecheckinghomomorphism_lemma_in_expression} S^*_3 S_3 + T^*_3 T_3 + T_3 T^*_3 + S^*_4 S_4 = 1.  \ee
 The lemma is derived by substituting the expressions of $ T^{*}_{3} T_3, ~ T_3 T^{*}_{3} $ and $ S^{*}_{4} S_4 $ from ( \ref{spherecheckinghomomorphism one_1} ), ( \ref{spherecheckinghomomorphism one_2} ) and Lemma \ref{spherecheckinghomomorphism two} in the equation ( \ref{spherecheckinghomomorphism_lemma_in_expression} ). \qed

\blmma

\label{spherecheckinghomomorphism four}

$ {( 1 + \mu^2 )}^{2} S_3 S^{*}_{3} = ( 1 + \mu^2 T_2 ) ( 1 + \mu^2 - \mu^4 ( 1 - T_2 ) ). $

\elmma

{\it Proof :} By unitarity of the matrix $ Z $, in particular equating the ( 1, 1 ) th  entry of $ Z Z^{*} $ to 1 we get $ S_3 S^{*}_{3} + \mu^2 ( 1 + \mu^2 ) T_3 T^{*}_{3} + \mu^2 S^{*}_{4} S_4 = 1 .$ Then the Lemma follows by using  ( \ref{spherecheckinghomomorphism one_2} ) and Lemma \ref{spherecheckinghomomorphism two}  in the above equation. \qed

\blmma

\label{spherecheckinghomomorphism five}

$ - S^{*}_{2} S_3 = ( \mu^2 + T_2 ) T_3. $

\elmma

{\it Proof :} By applying the adjoint and then multiplying by $ \mu^2 $ on ( \ref{spherehomomorphism B*B = A - A2 5} ) we have $ \mu^2 S^*_2 S_3 + \mu^2 S^*_4 S_2 = - \mu^2 T_3 ( \mu^2 + T_2 ) + \mu^2 ( 1 - T_2 ) T_3 .$ Subtracting this from ( \ref{spherehomomorphism B*B = A - A2 6} ) we have $ ( 1 - \mu^2 ) S^*_2 S_3 = - T_2 T_3 - \mu^2 T_3 T_2 + \mu^2 T_3 ( \mu^2 + T_2 ) - \mu^2 ( 1 - T_2 ) T_3 $ which implies $ - S^*_2 S_3 = ( \mu^2 + T_2 ) T_3  $ as $ \mu^2 \neq 1.$ \qed

\blmma

\label{spherecheckinghomomorphism six}

$ S_2 ( 1 - T_2 ) = \mu^2 ( 1 - T_2 ) S_2. $ 

\elmma 

{\it Proof :} Applying $ \kappa $  to  Lemma \ref{spherecheckinghomomorphism five} and then taking adjoint, we have
\be \label{spherecheckinghomomorphism six_1}  \mu^2 {( 1 + \mu^2 )}^{2} T^{*}_3 S_3 = - ( \mu^2 + T_2 ) S_{2}. \ee
Adding ( \ref{k on BB* = mu2 A - mu4 A2 B} ) and ( \ref{k on BB* = mu2 A - mu4 A2 AB} ) and then taking adjoint, we get ( by using $ \mu^2 \neq 1 $ ) 
\be \label{spherecheckinghomomorphism six_2}  \mu^2 {( 1 + \mu^2 )}^{2} T_3 S_4 = \mu^4 ( 1 - T_2 ) S_2. \ee
Moreover, ( \ref{k on B*B = A - A2 5} ) gives
$$ \mu^2 {( 1 + \mu^2 )}^{2} S_4 T_3 = - S_2 ( \mu^2 + T_2 ) + ( 1 - T_2 )S_2 - \mu^2 {( 1 + \mu^2 )}^{2} T^{*}_{3} S_3. $$
Using ( \ref{spherecheckinghomomorphism six_1} ), the right hand side of this equation turns out to be $ S_2 ( 1 - T_2 ).$

Thus,
\be \label{spherecheckinghomomorphism six_3}  {( 1 + \mu^2 )}^{2} S_4 T_3 = \mu^{-2} S_2 ( 1 - T_2 ). \ee
Again, application of adjoint to the equation ( \ref{k on BB* = mu2 A - mu4 A2 B} ) gives :
$$ \mu^2 {( 1 + \mu^2 )}^{2} S_3 T^{*}_3 = - \mu^2 {( 1 + \mu^2 )}^{2} T_3 S_4 - \mu^2 ( 1 + \mu^2 ) S_2 + \mu^4 ( 1 - T_2 )S_2 + \mu^4 S_2 ( 1 - T_2 ) .$$
Using ( \ref{spherecheckinghomomorphism six_2} ), we get
\be \label{spherecheckinghomomorphism six_4} {( 1 + \mu^2 )}^{2} S_3 T^{*}_3 = - S_2 ( 1 + \mu^2 T_2 ). \ee
Using ( \ref{spherecheckinghomomorphism six_1} ) - ( \ref{spherecheckinghomomorphism six_4} ) to the equation ( \ref{spherehomomorphism BA = mu2 AB A} ), we obtain :

$ S_2 - 2 S_2 T_2 - ( 1 + \mu^2 )^{-1} \mu^2 S_2 ( 1 + \mu^2 T_2 ) + \mu^{-2} ( 1 + \mu^2 )^{-1} S_2 ( 1 - T_2 ) = \mu^2 S_2 - 2 \mu^2 T_2 S_2 + ( 1 + \mu^2 )^{-1} \mu^6 ( 1 - T_2 )S_2 - ( 1 + \mu^2 )^{-1} ( \mu^2 + T_2 ) S_2. $

This gives

$ \mu^2 ( 1 + \mu^2 ) [ ( S_2 - S_2 T_2 ) - ( \mu^2 S_2 - \mu^2 T_2 S_2 ) ] - \mu^2 ( 1 + \mu^2 ) ( S_2 T_2 - \mu^2 T_2 S_2 ) - \mu^4 S_2 - \mu^6 S_2 T_2 + S_2 ( 1 - T_2 ) - \mu^8 ( S_2 - T_2 S_2 ) + \mu^4 S_2 + \mu^2 T_2 S_2 = 0 .$

Thus,
$ \mu^2 ( 1 + \mu^2 )[ S_2 ( 1 - T_2 ) - \mu^2 ( 1 - T_2 ) S_2 ] + S_2 ( 1 - T_2 ) - \mu^2 ( S_2 - T_2 S_2 ) + \mu^6 [ S_2 ( 1 - T_2 ) - \mu^2 ( 1 - T_2 ) S_2 ] - \mu^6 ( 1 - T_2 )S_2 + \mu^4 S_2 ( 1 - T_2 ) + \mu^2 ( S_2 ( 1 - T_2 ) - \mu^2 ( 1 - T_2 ) S_2 ) = 0 .$

On simplifying, $ ( \mu^6 + 2 \mu^4 + 2 \mu^2 + 1 )( S_2 ( 1 - T_2 ) - \mu^2 ( 1 - T_2 ) S_2 ) = 0 ,$ which proves the lemma as $ 0 < \mu < 1 .$ \qed

\blmma


\label{spherecheckinghomomorphism seven}


\be \label{spherecheckinghomomorphism seven_1} T_3 ( 1 - T_2 ) = \mu^2 ( 1- T_2 )T_3,  \ee
\be \label{spherecheckinghomomorphism seven_2} S_3 S^*_4 = \mu^4 S^*_4 S_3.  \ee

\elmma

{\it Proof :} The equation ( \ref{spherecheckinghomomorphism seven_1} ) follows by applying $ \kappa $ on Lemma \ref{spherecheckinghomomorphism six} and then taking adjoint.

We have $ S^*_4 S_3 = - T^2_3 $ from ( \ref{spherehomomorphism B*B = A - A2 7} ). On the other hand we have $ S_3 S^*_4 = - \mu^4 T^2_3 $ from ( \ref{spherehomomorphism BB* = mu2 A - mu4 A2 B2} ). Combining these two, we get ( \ref{spherecheckinghomomorphism seven_2} ). \qed

\blmma


\label{spherecheckinghomomorphism eight}

$ S_4 T_2 = T_2 S_4. $

\elmma

{\it Proof :} 
Subtracting ( \ref{k on BA = mu2 AB B*} ) from ( \ref{spherehomomorphism BA = mu2 AB B*} )  we get the required result. \qed

\blmma


\label{spherecheckinghomomorphism nine}

$ T_3 S_2 = S_2 T_3. $

\elmma

{\it Proof :} By applying adjoint on  ( \ref{k on BA = mu2 AB B} ) and then subtracting it from ( \ref{spherehomomorphism BA = mu2 AB B} ) we obtain $  S_2 T_3 - T_3 S_2  = 0 .$ 
\qed






\blmma


\label{spherecheckinghomomorphism ten}

$ S_3 ( 1 - T_2 ) = \mu^4 ( 1 - T_2 ) S_3. $ 

\elmma

{\it Proof :} By adding ( \ref{spherehomomorphism BA = mu2 AB B} ) with ( \ref{spherehomomorphism BA = mu2 AB AB} ) we obtain 

$$ S_3 ( 1 - T_2 ) + \mu^2 S_3 ( T_2 - 1 ) = \mu^2 ( \mu^2 - 1 ) T_3 S_2. $$ 
Thus, using $ \mu^2 \neq 1,$
\be \label{spherecheckinghomomorphism ten_1}  S_3 ( 1 - T_2 ) = - \mu^2 T_3 S_2.   \ee
 Moreover, by taking adjoint  of  ( \ref{k on BA = mu2 AB B} ), we obtain $ \mu^2 ( 1 - T_2 ) S_3 = \mu^2 S_2 T_3 - T_3 S_2 + S_3 ( 1 - T_2 ) .$
 
 Thus, $$ \mu^4 ( 1 - T_2 ) S_3 = \mu^4 S_2 T_3 - \mu^2 T_3 S_2 + \mu^2 S_3 ( 1 - T_2 ) .$$
 
 Hence, to prove the Lemma it suffices to prove: 
$$ S_3 ( 1 - T_2 ) = \mu^4 S_2 T_3 - \mu^2 T_3 S_2 + \mu^2 S_3 ( 1 - T_2 ). $$
After using $ T_3 S_2 = S_2 T_3 $ obtained from Lemma \ref{spherecheckinghomomorphism nine} we get this to be  the same as $ ( 1 - \mu^2 ) S_3 ( 1 - T_2 ) = \mu^2 ( \mu^2 - 1 ) T_3 S_2  .$ This is equivalent to $ S_3 ( 1 - T_2 ) = - \mu^2 T_3 S_2 $ ( as $ \mu^2 \neq 1 $ ) which follows from ( \ref{spherecheckinghomomorphism ten_1} ). \qed

\bppsn

\label{sphere_identification}

The map $ SO_{\mu}( 3 ) \rightarrow \clq $ sending $ M, L, G, N, C $ to $ - ( 1 + \mu^2 )^{-1} S_2,~ S_3,~ - \mu^{-1} S_4,~ ( 1 +\mu^2 )^{-1} ( 1 - T_2 ),~ \mu T_3 $ respectively is a CQG homomorphism.

\eppsn

{\it Proof :} It is enough to check that the map is $\ast$-homomorphic, since the coproducts on $SO_\mu(3)$ and $\clq$ are determined in terms of  
 the fundamental unitaries $Z^\prime$ and $Z$ respectively, and the map described in the statement of the proposition sends $(ij)$-th entry of $Z^\prime$ to the $(ij)$-th entry of $Z$ for all $(ij)$. 
 
 Now, it can easily be checked that the proof of the homomorphic property of the given map reduces to verification of the  relations on $ \clq $  as derived in Lemmas \ref{spherecheckinghomomorphism zero} - \ref{spherecheckinghomomorphism ten} along with the following equations :
\be \label{spherecheckinghomomorphism eleven} S_3 S_4 = \mu^4 S_4 S_3,   \ee
\be \label{spherecheckinghomomorphism twelve} S_3 S_2 = \mu^2 S_2 S_3, \ee
\be \label{spherecheckinghomomorphism thirteen} S_2 S_4 = \mu^2 S_4 S_2,  \ee
\be \label{spherecheckinghomomorphism fourteen}  S_3 S_4 = - \frac{\mu^2}{{( 1 + \mu^2 )}^2} S^2_2,  \ee
which follow from   Remark \ref{sphere_remark_k}, ( \ref{k on BA = mu2 AB B2} ), ( \ref{k on BA = mu2 AB B*2} ), ( \ref{k on BB* = mu2 A - mu4 A2 B2} ) respectively.   \qed

\bthm

\label{sphere_final_theorem}

We have the isomorphism:  $$  QISO^{+}_{R} ( \clo(S^{2}_{\mu,c}), ~ \clh, ~ D ) \cong SO_{\mu}( 3 ). $$

\ethm

{\it Proof :} $ SU_{\mu} ( 2 ) $ is an object in $ \widetilde{QISO^{+}_{R}}(D) $ as remarked before, and  thus  one gets a surjective 
 morphism  from $\widetilde{QISO^+_R}(D)$ to $SU_\mu(2)$ which clearly maps $QISO^+_R(D)$ onto $ SO_{\mu}( 3 ) $, identifying the latter as a quantum subgroup of $QISO^+_R(D)$. Let us denote the surjective map from $QISO^+_R(D)$ to $SO_\mu(3)$ by $\Pi$.
On the other hand,  Proposition \ref{sphere_identification} implies that $QISO^+_R(D)$ is a quantum subgroup of $SO_\mu(3)$, and the corresponding surjective CQG morphism from $SO_\mu(3)$ onto $QISO^+_R(D)$ is clearly seen to be the inverse of $\Pi$, thereby completing the proof. \qed 

\vspace{4mm}
\brmrk
Theorem \ref{sphere_final_theorem} shows that for a fixed $ \mu, $ the quantum isometry group $QISO^+_R(D)$ of   $ S^2_{\mu,c} $ does not depend on $c.$ This may appear somewhat surprising, but let us remark that in the classical situation ( that is for $ \mu = 1 $ ), $ c $ corresponds to the radius of the sphere and $S^2_{1,c}$ are isomorphic as $C^*$ algebras for all $c\geq0.$ We refer the reader to \cite{klimyk}, page 126, for the details regarding this. Although in the noncommutative case, that is, when $\mu \neq 1$, we do get non-isomorphic $C^*$-algebras $S^2_{\mu,c}$ for different choices of $c,$ one may still think that  the parameter $c$  in some sense determines the `radius' of the noncommutative sphere, and thus one should get the same (quantum) isometry group for different choices of $c.$ 

In view of the above, it seems impossible to `reconstruct' the quantum homogeneous spaces $S^2_{\mu,c}$ from the quantum isometry groups $SO_\mu(3)$. In this context, it may be mentioned that for $\mu \neq 1$, although all $S^2_{\mu,c}$ are quantum homogeneous spaces corresponding to $SO_\mu(3)$,  only $S^2_{\mu,0}$ arises as a quotient of $SO_\mu(3)$ by a quantum subgroup (see \cite{podles_subgroup} for more details). Thus, it is perhaps possible to  somehow `reconstruct'  $S^2_{\mu,0}$ from the quantum group $SO_\mu(3)$.
\ermrk

\subsection{ Existence of $ \widetilde{QISO^{+}}( D ) $ }

 For the above spectral triple, we have been unable to settle the issue of the existence of $ \widetilde{QISO^{+}}( D ) $ which is 
the universal object ( if it exists ) in the category $ {\bf{Q^{\prime}}} ( D ) $ mentioned in Section 1. Nevertheless, we shall show that if a universal object in ${\bf Q}^\prime(D)$  exists, then  $ QISO^{+}( D ) $ must coincide with $ SO_{\mu} ( 3 ).$

\blmma

\label{sphere_unresricted}

If $ \widetilde{QISO^{+}}( D ) $ exists, its induced action on $ S^{2}_{\mu,c} $, say $ \alpha_0 $, must preserve
 the state $h$ on the subspace  spanned by $ \{1,A,B,B^*,AB,AB^*,A^2,B^2,B^{*^2} \}.$

\elmma

{\it Proof :} Let $ \clw_0 =  \IC.1,$ $\clw_{\frac{1}{2}} = {\rm Span}  \{ 1, A, B, B^* \},$

 $ \clw_{\frac{3}{2}} ={\rm Span}
 \{ 1, A, B, B^*, AB, AB^*, A^2, B^2, B^{*2}  \}$. 

We note that the proof of Proposition \ref{sphere_Dabrowski_linearity_2} and the lemmas preceding it do not use the assumption that the action is $R$-twisted  volume 
 preserving, 
so the proof of Proposition \ref{sphere_Dabrowski_linearity_2} goes through verbatim implying that 
$ \alpha_0 $ keeps invariant the subspace spanned by $ \{ 1, A, B, B^* \} $ and hence 
it preserves  $ \clw_{\frac{3}{2}}$ as well.   Let $\clw_{\frac{3}{2}}=\clw_{\frac{1}{2}} \oplus \clw^\prime$ be the orthogonal decomposition with respect to the Haar 
 state (say $h_0$) of $QISO^+(D)$. Since
 $ SO_{\mu} ( 3 ) $ is a sub-object of $QISO^+(D)$,
 there is a CQG morphism $ \pi $ from $ QISO^{+}( D ) $ onto $ SO_{\mu} ( 3 ) $ satisfying $ ( {\rm id} \otimes \pi ) \alpha_0 = \Delta$, where 
 $\Delta$ is the $SO_\mu(3)$ action on $S^2_{\mu,c}$. It follows from this that any $QISO^+(D)$-invariant subspace 
(in particular $\clw^\prime)$ is also $SO_\mu(3)$-invariant. On the other hand, it is easily seen that on $\clw_{\frac{3}{2}}$, the $SO_\mu(3)$-action decomposes 
  as $\clw_{\frac{1}{2}} \oplus \clw^{\prime \prime}$, (orthogonality with respect to $h$, the Haar state of $SO_\mu(3))$), where $\clw^{\prime \prime}$ is a five 
 dimensional irreducible subspace.

We claim that $\clw^\prime=\clw^{\prime \prime}$, which will prove that the $QISO^+(D)$-action $\alpha_0$ has the same $h$-orthogonal decomposition 
 as the $SO_\mu(3)$-action on $\clw_{\frac{3}{2}}$, so preserves $\IC.1$ and its $h$-orthogonal complements. This will prove that $\alpha_0$ preserves the Haar 
 state $h$ on $\clw_{\frac{3}{2}}$. 

We now prove the claim. Observe that $\clv:=\clw^\prime \bigcap \clw^{\prime \prime}$ is invariant under the $SO_\mu(3)$-action but due to the irreducibility of $ \Delta $ on the vector space $ \clw^{\prime} $ or $ \clw^{\prime \prime} ,$ it has to be zero or  $ \clw^{\prime} = ~ \clw^{\prime \prime} .$ Now,  $ {\rm dim} ( \clv ) = 0 $ implies $ {\rm dim} ( \clw^{\prime} ) + {\rm dim} ( \clw^{\prime \prime} ) = 5 + 5 > 9 = {\rm dim} ( \clw_{\frac{3}{2}} ) $ which is a contradiction unless $ \clw^{\prime} = ~ \clw^{\prime \prime} .$

\qed 

 \bthm
\label{existence_issue}
If $ \widetilde{QISO^{+}}( D ) $ exists, then we must have that  $ QISO^{+} ( D ) \cong SO_{\mu} ( 3 ). $
\ethm

{\it Proof :} In the proof of Lemma \ref{sphere_unresricted}, it was noted that Proposition \ref{sphere_Dabrowski_linearity_2} follows under the assumption of the present theorem. 
 To complete the proof of the theorem, we just need to observe  that the other lemmas used  for proving Theorem \ref{sphere_final_theorem} require the conclusion of Lemma \ref{sphere_unresricted} as the only extra ingredient. \qed

Let us conclude the article with brief explanation of   the technical difficulties regarding the issue of existence of $\widetilde{QISO^{+}}( D )$ . Let $ R^{\prime} $ be a positive, invertible operator commuting with $ D $ such that $ \tau_{R^{\prime}} \neq \tau_{R} $ and let $ \phi^{\prime} $ denote the action of $ QISO^{+}_{R^{\prime}}(D) $ on $ S^2_{\mu,c}. $ The problem of existence of  $ \widetilde{QISO^+} ( D ) $ is closely related to   the question whether it is possible to identify  $ QISO^{+}_{R^{\prime}} ( D ) $ as a quantum subgroup  of $ SO_{\mu} ( 3 ) $ for a general $ R^{\prime}. $  By Theorem \ref{existence_issue}, a negative answer of this question will prove that ${\bf{Q^{\prime}}} ( D )$ does not have  a universal object.

 Now, as has been noted in Remark \ref{linearity_without_r},  $ \phi^{\prime} $ is linear, that is, it keeps the span of $ \{ 1, A, B, B^* \} $ invariant and hence  it us given by an expression similar to equations ( \ref{alphaA} ) and ( \ref{alphaB} )  with $ T_i, S_i $ replaced by some $ T^{\prime}_i, S^{\prime}_i $ which generate   $ QISO^{+}_{R^{\prime}} ( D ) $ as a $C^*$-algebra. We can in principle write down all the relations satisfied by these generators, proceeding as in the Subsection 3.3. These relations will be analogous to  equations ( \ref{spherehomomorphism B*B = A - A2 1} )-(  \ref{k on BB* = mu2 A - mu4 A2 AB} ), and in fact, the relations which make use of the homomorphism property only remain unchanged. However, the ones making use of the fact that $\phi^\prime$ preserves $\tau_{R^\prime}$ will change, since $\tau_{R^\prime}$ is in general different from $\tau_R$. 
  In particular, the expression of the antipode will change, which will affect all the relations starting from ( \ref {k on B*B = A - A2 2} ). We need to have a  deeper and systematic understanding of the relations satisfied by $T^\prime_i, S^\prime_i$    for a general $ R^{\prime}, $ and possibly study their representations 
    in concrete Hilbert spaces, to decide whether $QISO^+_{R^\prime}(D)$ is a quantum subgroup of $SO_\mu(3)$ or not.  We are not yet able to do this.
    
    Moreover, even if $\widetilde{QISO^+}(D)$ exists,  although we can  identify $QISO^+(D)$ with the well-known quantum group 
     $SO_\mu(3)$,  it is not so easy to explicitly compute  $\widetilde{QISO^+}(D)$. If $U$ denotes the unitary representation corresponding to $\widetilde{QISO^+}(D)$, the fact that $U$ commutes with $D$ implies that $U$ must preserve each of the two-dimensional eigenspaces $ {\rm span} \{ v^{l}_{m, \frac{1}{2}} + v^{l}_{m, - \frac{1}{2}} : m = \pm \frac{1}{2} \} $ and $ {\rm span} \{ v^{l}_{m, \frac{1}{2}} - v^{l}_{m, - \frac{1}{2}} : m = \pm \frac{1}{2} \} $ of $D$. Suppose that $( q_{ij})_{i,j=1,2}$ and $( r_{ij} )_{i,j=1,2}$ are the matrices (with entries in $\widetilde{QISO^+}(D)$ ) of $U$ corresponding to these two spaces respectively. Then it is clear that as a $C^*$ algebra $\widetilde{QISO^+}(D)$ will be generated by $q_{ij}, r_{ij}$'s as well as the generators $T_i, S_i$ of $SO_\mu(3)$. However, the mutual relations among these generating elements have to be determined from the fact that $U$ preserves each of the eigenspaces of $D$. In principle one gets infinitely many such relations which are quite complicated and it is not clear how to simplify them.

Jyotishman Bhowmick : Stat-Math Unit, Indian Statistical Institute, 203, B. T. Road, Kolkata 700 108
E-mail address: jyotishmanb@gmail.com

Debashish Goswami: Stat-Math Unit, Indian Statistical Institute, 203, B. T. Road, Kolkata 700 108
E-mail address: goswamid@isical.ac.in

\end{document}